\documentclass[12pt]{amsart}

\usepackage{graphicx, verbatim,amsmath,amssymb, float}
\usepackage{hyperref}
\hypersetup{
    colorlinks=true,
    linkcolor=blue,
    filecolor=magenta,      
    urlcolor=cyan,
}

\def\vs#1 {\vskip#1truein}
\def\hs#1 {\hskip#1truein}

\newcommand{\be}{\begin{equation}}
\newcommand{\ee}{\end{equation}}
\newcommand{\R}{{\mathbb R}}

\newtheorem*{thm*}{Theorem}

\everymath{\displaystyle}

\theoremstyle{definition}

\title{Nucleation during phase transitions in random networks}
\date{\today}

\author{
Joe Neeman
\and Charles Radin
\and Lorenzo Sadun
}
\address{Joe Neeman\\Department of Mathematics\\The University of
  Texas at Austin\\ Austin, TX 78712} \email{joeneeman@gmail.com}
\address{Charles Radin\\Department of Mathematics\\The University of
  Texas at Austin\\ Austin, TX 78712} \email{radin@math.utexas.edu}
\address{Lorenzo Sadun\\Department of Mathematics\\The University of
  Texas at Austin\\ Austin, TX 78712} \email{sadun@math.utexas.edu}

\begin{document}

\begin{abstract}
We analyze the 3-parameter family of random networks that are uniform
on networks with fixed number of edges, triangles, and nodes (between
33 and 66). We find precursors of phase transitions which are
known to be present in the asymptotic node regime as the edge and
triangle numbers are varied, and focus on one of
the discontinuous ones. By use of a natural edge flip dynamics we determine
nucleation barriers as a random network crosses the transition, in
analogy to the process a material undergoes when frozen or melted, and characterize
some of the stochastic properties of the network nucleation.
\end{abstract}

\maketitle

\section{Introduction}\label{introduction}

When a fluid is frozen to a
crystal a disordered array of molecules has to rearrange 
into an ordered array.
In practice molecules have difficulty doing this; careful experiments
can easily cool a liquid into a disordered supercooled liquid, well
below the freezing temperature, before transforming it into a
crystalline solid~\cite{Ca}.  The manner in which the difficulty is eventually
overcome is traditionally modelled by Classical Nucleation Theory,
CNT, some of whose predictions are off by many orders of magnitude~\cite{Ca}.

There is no satisfactory solvable toy model that exhibits a
fluid/solid phase transition - this is a famous unsolved problem~\cite{Br,Uh}
- so we will use a mathematical model of random networks, which has
adjustable parameters and exhibits phase transitions with the
parameters treated as analogues of thermodynamic parameters. The model
is unphysical but may still be useful to determine precisely how a
dynamically driven phase transition between incompatible structures is
actually obtained.

Our model is a mean-field model, namely the random `edge/triangle' network model
popularized by Strauss in 1986~\cite{St}. We first give some
background to explain the meaning of the transitions found in these
models. We then sketch our results on dynamically driving such a
system between phases, exhibiting a form of nucleation as 
incompatible structures are bridged.

\smallskip

We begin with some notation.  Given $n$ abstract labelled nodes, we
denote by $E(G)$ and $T(G)$ the number of edges and triangles in the
simple network $G$. (Simple means that two nodes can have at most one
edge connecting them and there are no loops connecting a node with
itself.)  Strauss put on the set of all such networks the 2-parameter
family of unnormalized probability mass functions: \be \exp[s_1 e(G) +
  s_2 t(G)],
\label{Strauss}
\ee
where 
\be
e(G)={E(G)\over {n \choose 2}}
\hbox{ and } t(G)={T(G)\over {n \choose 3}}
\label{defs}
\ee
 are the edge and
triangle densities in $G$, with value 1 in the complete network. The
model parameters $s_1, s_2$ are analogues of $\beta \mu$ and $\beta$
in a grand canonical distribution of statistical mechanics, where
$\beta$ is inverse temperature and $\mu$ is chemical potential.

A network $G$ is described by its adjacency matrix $M(G)$. This is 
a symmetric $n\times n$ matrix whose $(i,j)$ entry is 1 if nodes $i$ and $j$ 
are connected by an edge, and zero otherwise. The statistical properties of 
the matrix are independent of the (arbitrary) ordering of the nodes. 
Permutation-invariant data such as the eigenvalues of $M(G)$ are 
particularly informative. 

As in statistical mechanics, phases in parametric models of random
networks require an infinite-system-size limit formalism. Instead of looking
at a symmetric 0--1 matrix $M(G)$, one considers a symmetric function 
$g$ on the unit square $[0,1]^2$ taking values in the interval $[0,1]$. 
(The points in $[0,1]$ are thought of as representing infinitely many
nodes.)
Such a {\em graphon} $g$ is the kernel of an integral operator 
$M(g)$ acting on $L^2([0,1])$. 
In an important
series of papers~\cite{BCL,BCLSV,LS1,LS2,LS3}, Borgs et al 
used a certain metric on the space
of networks to construct a completion, called the space of (reduced)
graphons, with many useful properties. Used together with an important
large deviation result~\cite{CV} of Chatterjee/Varadhan, one may study
convergence of distributions to infinite size through
optimal-free-energy graphons $g_{s_1,s_2}$. See for instance~\cite{CD}.

These network models are mean-field in the sense that there is no
space cutoff limiting the direct influence of variables on one another. For
traditional statistical mechanics models with short range forces there is an
equivalence of ensembles, connected by Legendre transforms.  However,
as is common in mean-field models, Legendre transforms are not invertible 
for these graphon models and lose
information contained in the microcanonical ensemble. This is not
necessarily a problem when studying a network with a small number of nodes, as
is common in social science modelling~\cite{Ne}. However, when studying
asymptotically large networks to study phase transitions the loss would be a
problem so we
must use the microcanonical version of the
model~\cite{KRRS1, KRRS2}. That is, we must specify the parameters $e(G)$ and 
$t(G)$ and consider the uniform distribution on networks with the fixed
values of $e$ and $t$. We then study convergence to infinite
size through optimal-entropy graphons, $g_{e,t}$.  

For fixed values of $(e,t)$, the graphon formalism allows one to
take a limit in the number of nodes, $n\to \infty$, with `most'
constrained networks $G$ converging to an entropy-optimal graphon
$g_{(e,t)}$, in the sense that their (suitably normalized) adjacency matrices
$M(G)$, viewed as operators on $\R^n$, converge in a weak
sense to the compact linear integral operator $M(g_{(e,t)})$ on
$L^2([0,1]^2)$, with the (suitably normalized) 
ordered sets of eigenvalues of the $M(G)$'s
converging to those of $M(g_{(e,t)})$; see section 11.6 in \cite{Lo}.
In~\cite{KRRS1,KRRS2,KRRS3,Ko,RS1,RS2,RRS1,RRS2,RRS3} this and similar
models have been studied in depth; we describe some results next. 

The optimal-entropy graphons have been shown
to have a very simple form. 
If we break the interval $[0,1]$ into finitely many sub-intervals,
called ``podes'',
and thereby break the square $[0,1]^2$ into finitely many rectangular
blocks, then the optimal graphon is constant on each block. As noted
above, the
interval $[0,1]$ is thought of as representing infinitely many nodes,
so we are envisioning a decomposition of the set of nodes into
a finite number of collections, called podes, of `equivalent' nodes, 
equivalent in
the sense that two equivalent nodes have the same probability of
connection with any other node.
It only takes a few parameters
to describe such a ``multipodal''  graphon, namely the sizes of the
sub-intervals and the value of $g$ on each block;
this vast reduction in degrees of freedom is
chiefly responsible for the proliferation of results in these
models. 

By a 
{\em phase} in our network model we mean an open connected region in 
the space
of pairs
$\{(e, t)\}$, at each point of which there is a unique optimal-entropy
graphon $g_{(e,t)}$, whose parameters are {\em smooth} functions of $(e,t)$. 
At a phase boundary there may be two or more optimal-entropy graphons,
typically with qualitatively different structures, 
in which case we say the phase transition is {\em
  discontinuous}. Alternatively there
may be a single optimal-entropy graphon, in which case we say the
phase transition is {\em continuous}. These are in some sense analogous to
first-order and second-order phase transitions in statistical mechanics, 
but the analogy can only be taken so far, especially since this
is a mean-field model. (In particular, the expansion
of the entropy as a power series in $e$ and $t$ sometimes involves fractional
powers.)

\begin{figure}[H]
\center{\includegraphics[width=5in]{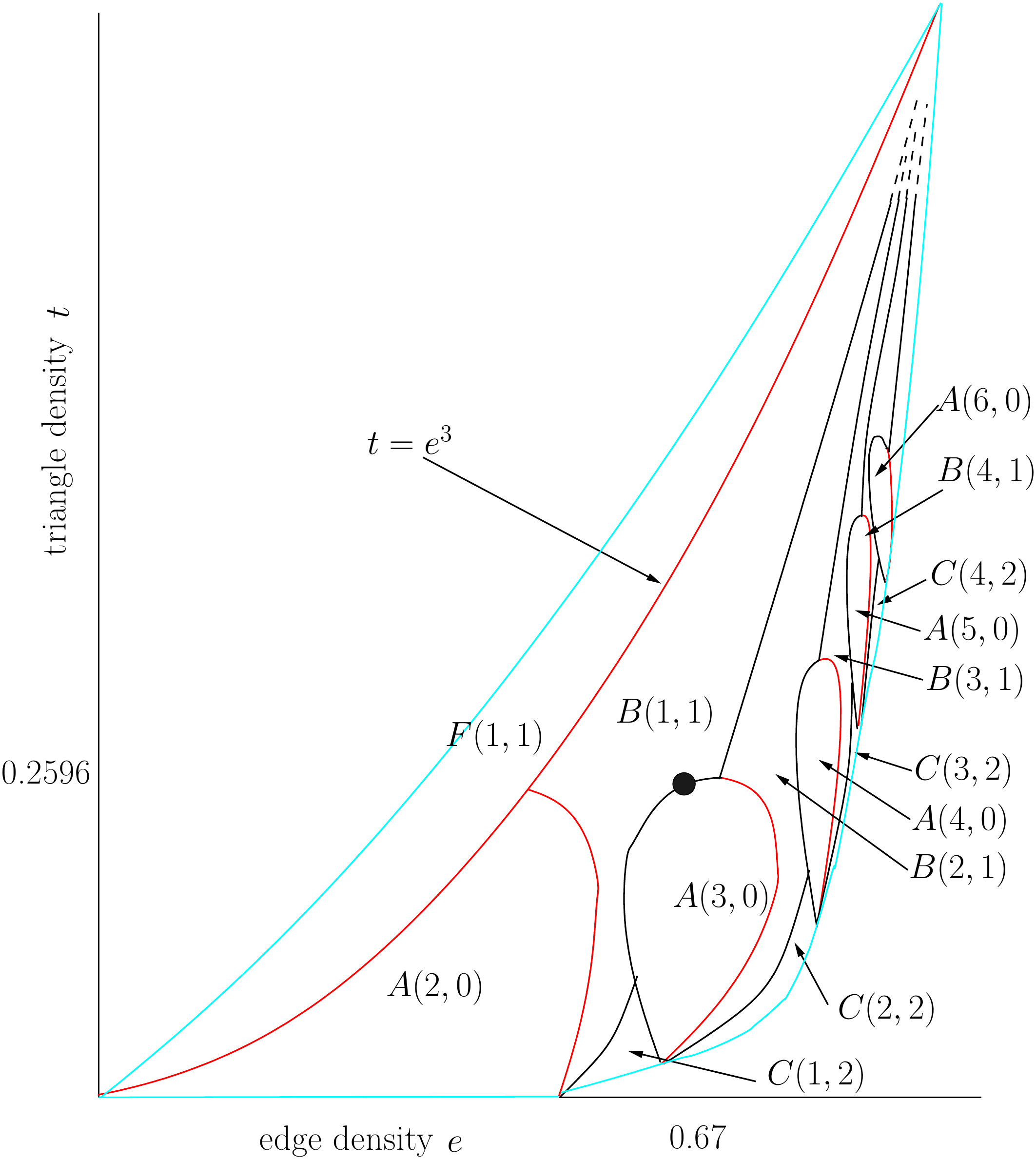}}
\vskip .1truein
\caption{{Schematic} sketch of $14$ of the phases in the edge/triangle
  model. Continuous phase transitions are shown in red, and
  discontinuous phase transitions in black. The boundary of the phase
  space is in blue. The transition point studied in this paper is indicated
  by a black dot.}
\label{phase-diagram}
\end{figure}

\begin{figure}[H]
\center{\includegraphics[width=5in]{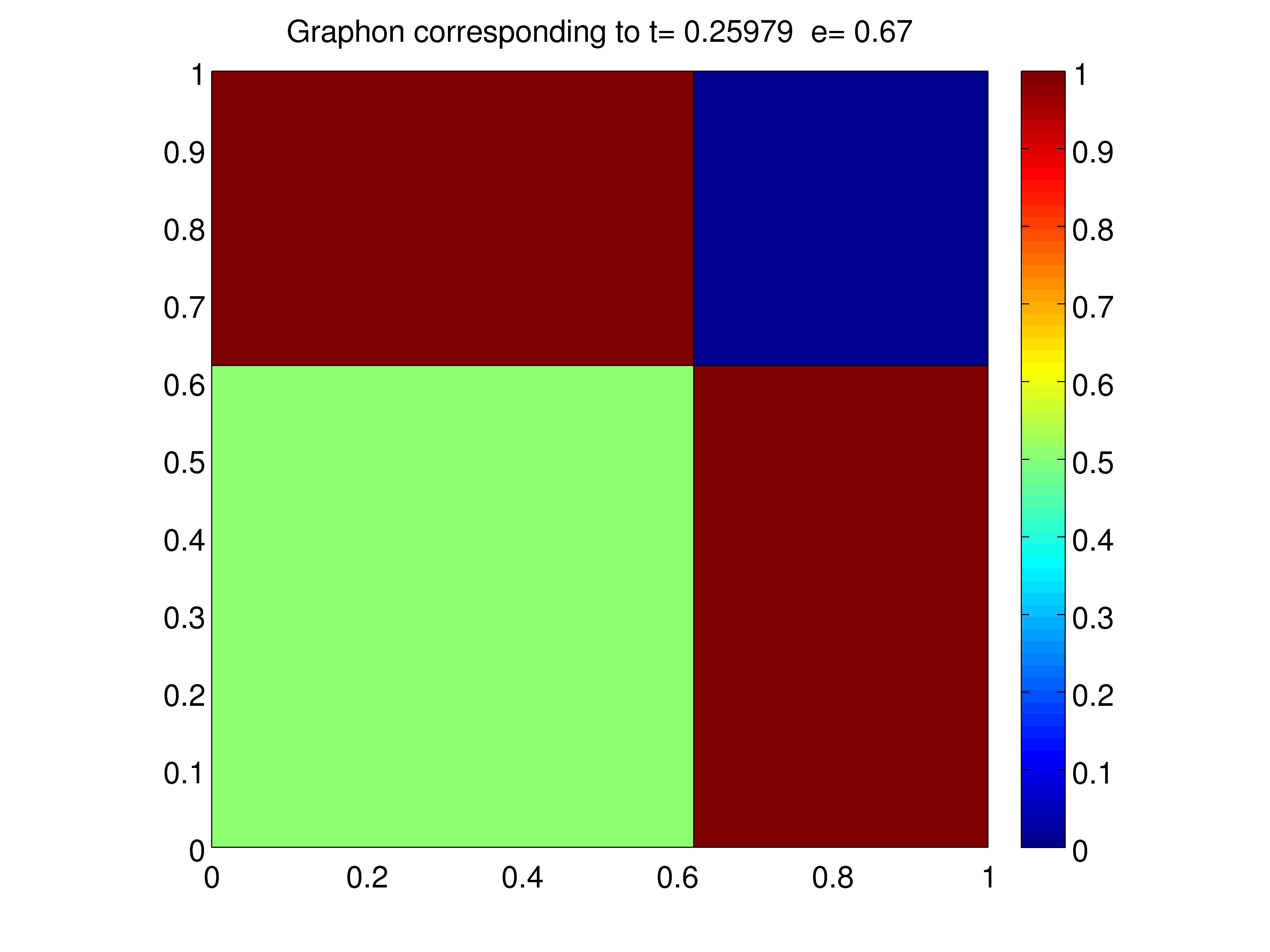}}

\center{\includegraphics[width=5in]{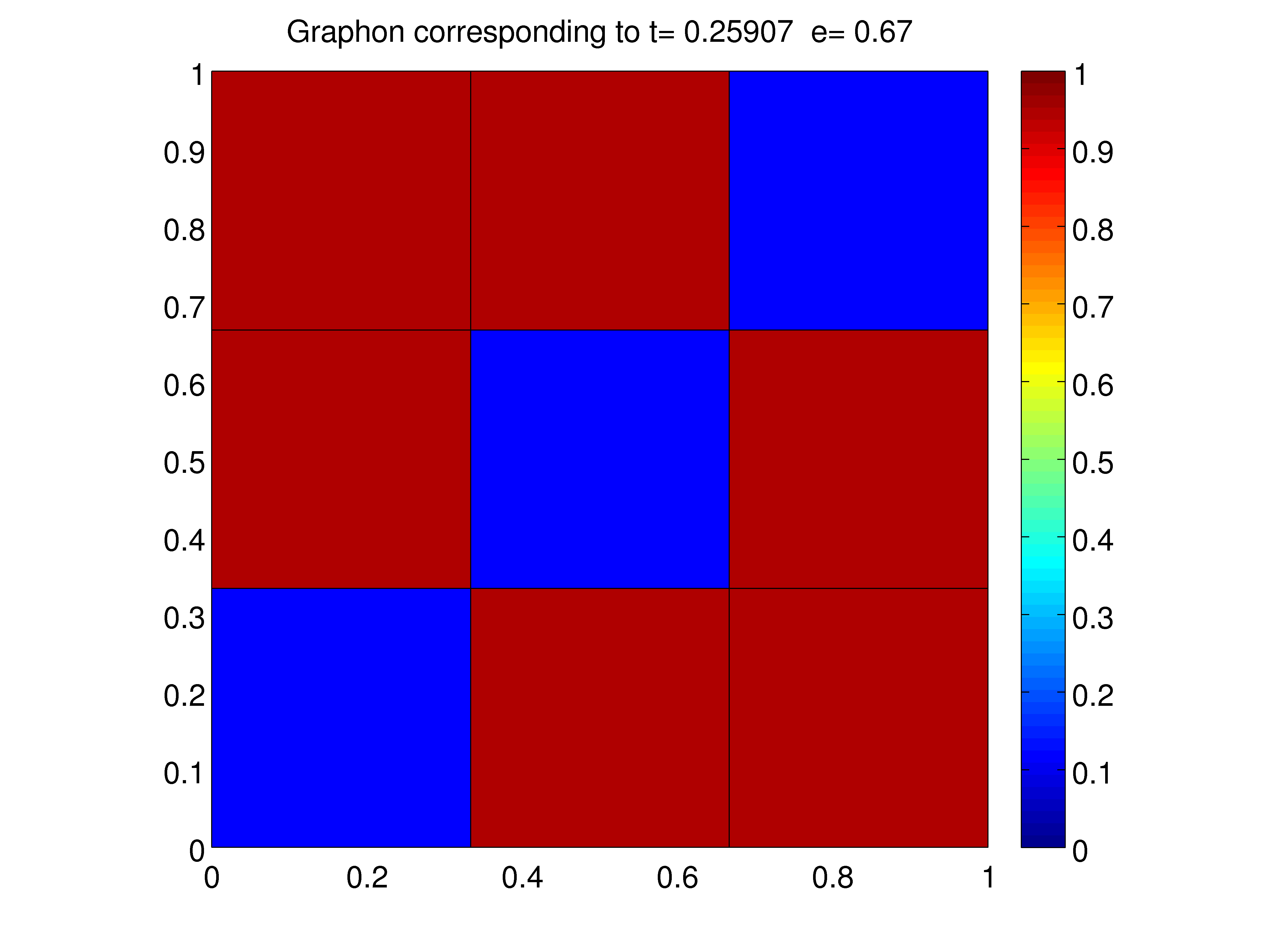}}
\caption{Transition between phases $B(1,1)$ and $A(3,0)$ at edge density
$0.670$.
Triangle densities of the two graphons: 0.25979, 0.25907.
Triangle density at transition: 0.25963868.}

\label{optimalgraphons}
\end{figure}

In Figure~\ref{phase-diagram} we sketch the phase diagram of the 
microcanonical edge/triangle model~\cite{KRRS1}. 
This model exhibits three infinite families of phases labelled
$A(m,0)$, $B(m,1)$ and $C(m,2)$ as well as a singleton $F(1,1)$. 
An $A(m,0)$ graphon (that is, the optimizing graphon in the $A(m,0)$ phase)
has $m$ podes and is invariant under permutation of
these indistinguishable podes; see Figure~\ref{optimalgraphons}.
A $B(m,1)$ graphon has $m$ indistinguishable
podes and one pode that is different. A $C(m,2)$ graphon has $m$ 
indistinguishable podes and two additional 
podes that are different from the $m$ podes but are indistinguishable 
from each other. The differences between these phases are reflected in
the different eigenvectors of the adjacency matrices of generic networks in
each phase. 

We will be concentrating on the transition between the 
$B(1,1)$ and $A(3,0)$ phases as we go through the point 
$(e,t)= (0.67, 0.2596)$ along constant edge density
$0.67$. This transition is particularly interesting because it 
has been proven to be discontinuous; the optimizing graphons for $t$ slightly
above $0.2596$ are very different from the optimizing graphons 
for $t$ slightly below
$0.2596$.  (See Figure \ref{optimalgraphons}.) 

In switching from an
arrangement described by one graphon to an arrangement described by the 
other, nodes need to rearrange their connections in a process analogous
to nucleation in a liquid/crystal transition. Patterns are broken, patterns are formed, and once 
sufficiently many nodes follow the new pattern, the rest follow suit. 

Our study of this transition involved two major
challenges. First, we needed to work with finite networks, yet the very
notion of phases requires taking infinite-size limits. Although
the graphon formalism ensures that systems of large enough size are
modelled well by graphons, there is no mechanism in the formalism to 
determine how large is `large enough'. 
Fortunately, we were able to see
unambiguous evidence of the emergence of the $B(1,1)/A(3,0)$
transition in systems with as few as 30 nodes, well below the 100-or-so 
nodes that our computers can comfortably handle; we describe this in
Section~\ref{equilibrium}. The fact that such transitions can be seen
in systems of such low size is significant both for the computer
analysis of network transitions and for the modelling of
scientific data of such size, and is the {\bf first major result} of this paper.

Second, we needed to introduce a random dynamics to play the role of
cooling or heating a material 
from an equilibrium state at one set of thermodynamic parameters to
an equilibrium state at another set.
Our dynamics respects the constraint on the edge 
density while introducing enough randomness to drive the triangle
density across a transition.

The solutions we found are described in Section~\ref{nucleation}, both for 
going from $B(1,1)$ to $A(3,0)$ and for going from $A(3,0)$ to $B(1,1)$.
Along the way we produced videos, available at 

\href{https://www.youtube.com/playlist?list=PLZcI2rZdDGQqx3WEY6BXoJXcxqQXF5ZzQ}{https://www.youtube.com/playlist?list=PLZcI2rZdDGQqx3WEY6BXoJXcxqQXF5ZzQ}

\noindent that literally show nodes rearranging into the pode structure
appropriate to prescribed constraints. (The rearrangement is not in physical
space, which plays no role in this mean-field model. Rather, the podal structure
and its changes are seen in the eigenvectors of the adjacency matrix,
as noted above.)

Understanding this process is the {\bf second major result} of this paper. 
In each direction, the evolution from one phase structure
to the other occurs in three distinct stages. The first stage is very quick
and brings the density $t$ of triangles to its desired value. 
The second stage is
a factor of 10--100 slower and results in a structure that can be viewed as 
``new phase plus defects'', e.g. with podes that are slightly too big or too
small. The third phase is a factor of 1000--10,000 slower than the first and
achieves true equilibrium.  

Finally, the time required for 
crystal nucleation in physical systems is random, as are other features
of the process. We have begun an statistical analysis of the time
required for each of the three stages of network nucleation (both for
$A(3,0) \to B(1,1)$ and for $B(1,1) \to A(3,0)$), 
with details in Section~\ref{nucleation}.

\section{Equilibrium in finite systems}\label{equilibrium}

Our goal is to understand how a network changes from the structure of
one phase to that of another, under an appropriate dynamics. 
Before we can do this, we must
understand the structure of the two phases. We begin with the entropy-maximizing
graphons that describe infinite (or sufficiently large) random networks, and then consider finite-size
effects. For more about these graphons and the precise transition between 
them, see~\cite{KRRS1}.

In the A(3,0) phase, an optimizing graphon $g_{e,t}$ is shown in the bottom of 
Figure~\ref{optimalgraphons}. If we denote the value of the graphon along the diagonal by $a$ and
the off-diagonal value by $b$, then this describes an infinite network with 
the following properties. 
\begin{enumerate}
\item The nodes group into three classes or ``podes'' of equal
size, which we might imagine as red, blue, and green nodes. 
\item The edges within each pode, e.g. connecting two red nodes, appear
with probability $a$. 
\item The edges between different podes appear  with probability $b>a$. 
\item The corresponding integral operator $M(g_{e,t})$ has rank three. Its
  nonzero eigenvalues are $(a+2b)/3$ with
  multiplicity 1 and $(a-b)/3$ with multiplicity 2.
\item The positive eigenvalue $(a+2b)/3$ is the overall edge density, 
and the corresponding eigenfunction is constant. 
\item The triangle density is 
$$ t = \frac19 a^3 + \frac23 ab^2 + \frac29 b^3 = e^3 + 2 \left(\frac{a-b}{3}
\right)^2,$$
which is the trace of the cube of the integral operator. Given $e$ and $t$,
it is easy to compute $a$ and $b$, the eigenvalues and eigenfunctions of
$M(g_{e,t})$, and the densities of all subnetworks.
\item The eigenspace corresponding to the negative eigenvalue
$(a-b)/3$ consists of functions that are constant on each pode and 
integrate to zero. If we pick an orthonormal basis $\{\xi_1(x), \xi_2(x)\}$ 
for this eigenspace and pick a random point $x \in [0,1]$, 
the ordered pair $(\xi_1(x), \xi_2(x))$ can lie at three possible points 
in $\R^2$, and these points form the vertices of an equilateral triangle
centered at the origin. 
A change of basis only serves to rotate the triangle. 
\end{enumerate}

If we pick a network $G$ with $n$ nodes whose adjacency matrix $M(G)$ is 
close to this graphon (which necessarily requires $n$ large),
then  $M(G)$ will 
have two large negative eigenvalues of approximate size $n(a-b)/3$ and one 
large positive eigenvalue of approximate size $n(a+2b)/3$. The eigenfunctions
for the negative eigenvalues will be approximately constant on each pode,
so that if we plot the values $\{(\xi_1(v_i),\xi_2(v_i))\}$ for the $n$ nodes
$\{v_i\}$, we will get (with high probability) 
three clusters of roughly equal size,
centered at the vertices of an equilateral triangle. Moreover, the second
largest negative eigenvalue $\lambda_2$ will be nearly as large as the largest
negative eigenvalue $\lambda_1$, while all other negative eigenvalues will
be of order $\sqrt{n}$. 

We now turn to the $B(1,1)$ phase, with an optimal graphon shown at
the top of Figure~\ref{optimalgraphons}. The corresponding integral
operator has rank two, and in particular only has one negative
eigenvalue, whose eigenfunction is constant on each
pode. If we pick a large random network that is close to this
structure, then the most negative eigenvalue $\lambda_1$ 
of the adjacency matrix will be $O(n)$, and the values of
$\xi_1(v)$ will distinguish the two podes, but the second most 
negative eigenvalue $\lambda_2$ will only
be $O(\sqrt{n})$. $\xi_2$ will be essentially zero on the smaller pode
and random on the larger pode. A scatter plot of $\{(\xi_1(v_i),
\xi_2(v_i))\}$ will then give a tight cluster of points corresponding
to the small pode and a vertical bar corresponding to
the large pode.

The upshot is that both the size of $\lambda_2$ and the distribution
of $(\xi_1(v), \xi_2(v))$ are strong indicators of the phase we are
in. Another indicator is the size of the podes
themselves, as indicated by the size of the clusters in the 
$(\xi_1,\xi_2)$ plot. In the $A(3,0)$ phase, 
there are three podes of equal size. In
$B(1,1)$, there are two podes with a roughly 60-40 split between the
larger and the smaller pode.

\begin{figure}[htbp]
\center{\includegraphics[width=4in]{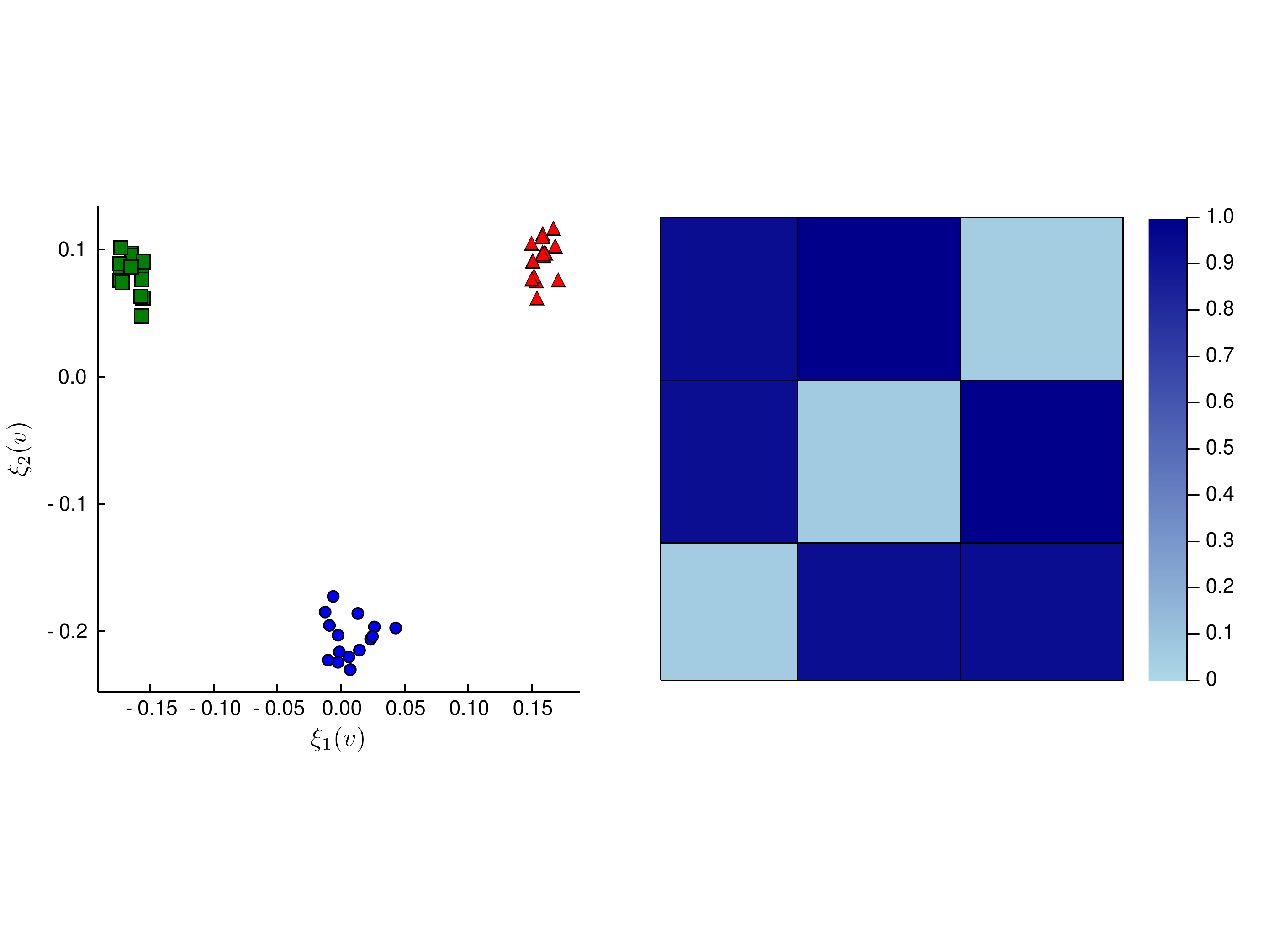}}
\caption{$(\xi_1(v),\xi_2(v))$ plot and best-fit tripodal graphon for a typical graph
    with $54$ nodes and triangle density $t=0.24$.  Each mark on the left
    corresponds to a node in the graph. Using spatial clustering on the
    left-hand plot, these nodes were divided into three podes (denoted by
    different colors and shapes in the left-hand plot).
    The plot on the right shows the empirical graphon corresponding to this division into
    three: the division of the axes represents the three podes, and the color of a region
    represents the density of edges between the corresponding podes.
}
\label{xiplot-a30}
\end{figure}

\begin{figure}[htbp]
\center{\includegraphics[width=4in]{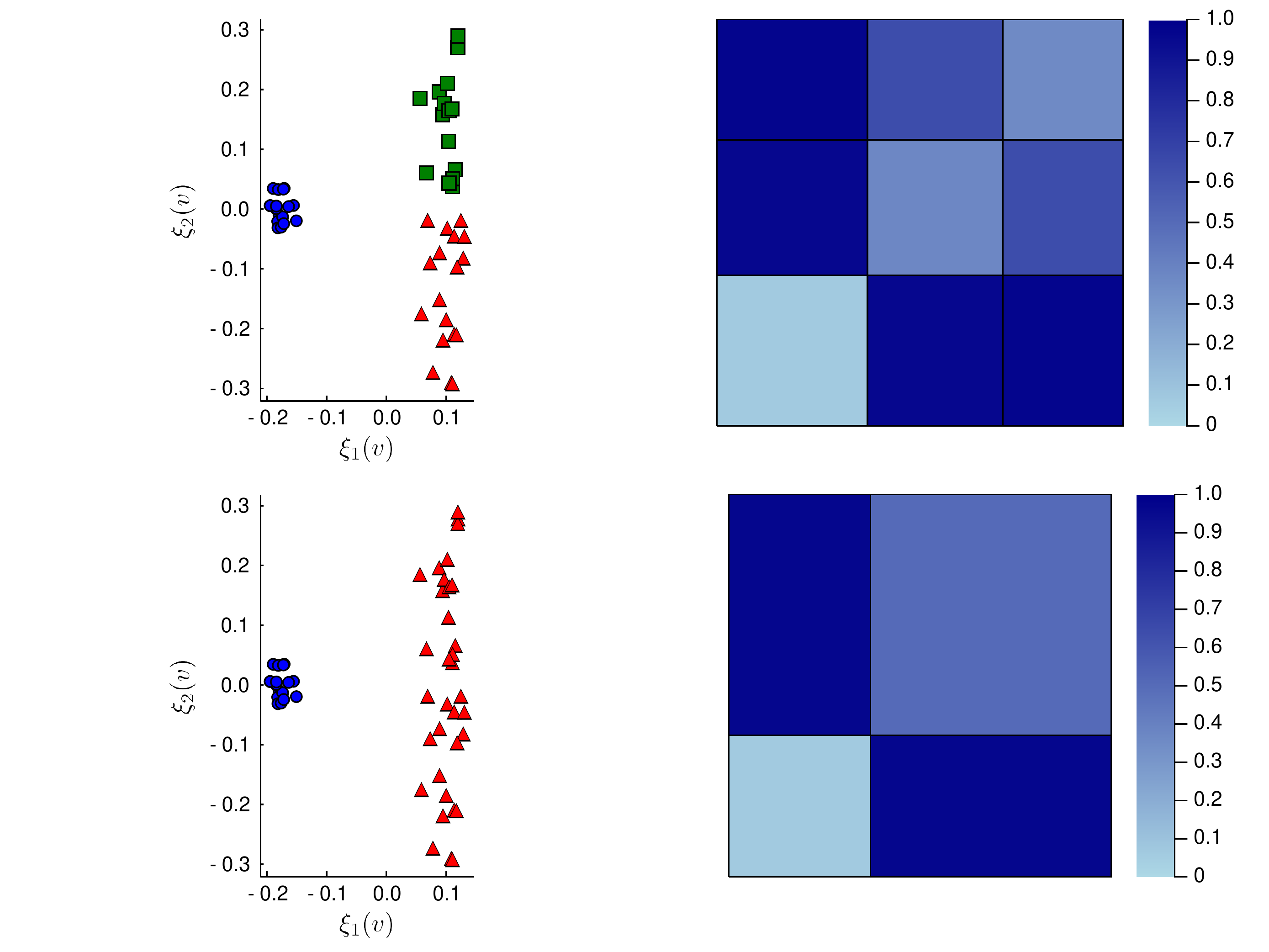}}
\caption{$(\xi_1(v),\xi_2(v))$ plot and best fit tri- and bipodal graphons for a typical
    graph with $54$ nodes and triangle density $t=0.26$. See Figure~\ref{xiplot-a30} for
    a thorough description.}
\label{xiplot-b11}
\end{figure}

The left side of Figure \ref{xiplot-a30} shows the distribution of
$(\xi_1(v), \xi_2(v))$ for a typical graph with $n=54$, $e=0.67$ and
$t=0.24$. Grouping nodes, we then generate a tripodal graphon by
simply counting the number of edges that blue-blue, blue-green,
etc. Figure \ref{xiplot-b11} is similar, only for a sample graph at
$t=0.26$, which is in the B(1,1) phase. In the plot of 
Figure~\ref{xiplot-b11}, we
attempt to segregate the nodes into three podes, as before, but there
is no clear dividing line between the reddish triangles and
squares. We prefer to treat such
graphs as in the second plot, where the entire vertical bar is
considered a single pode.

Of course, the uniform ensemble of networks of a fixed node number $n$
and fixed numbers of edges and triangles are only proven to look like
the optimal graphons in the $n \to \infty$ limit. We need a procedure
to generate the ensemble for finite $n$.  To do this, we implement a
dynamics on the space of
networks with given $e$, described in Subsection~\ref{sampling}, 
which has as limit set the subset of those
with given $t$. (The dynamics then plays the role of an 
adjustable heat bath.) We
then run the dynamics `sufficiently long' to obtain samples of this
limit set.  We next examine the histogram of values of $\lambda_2$ we
obtained this way.

\begin{figure}[htbp]
\center{\includegraphics[width=5in]{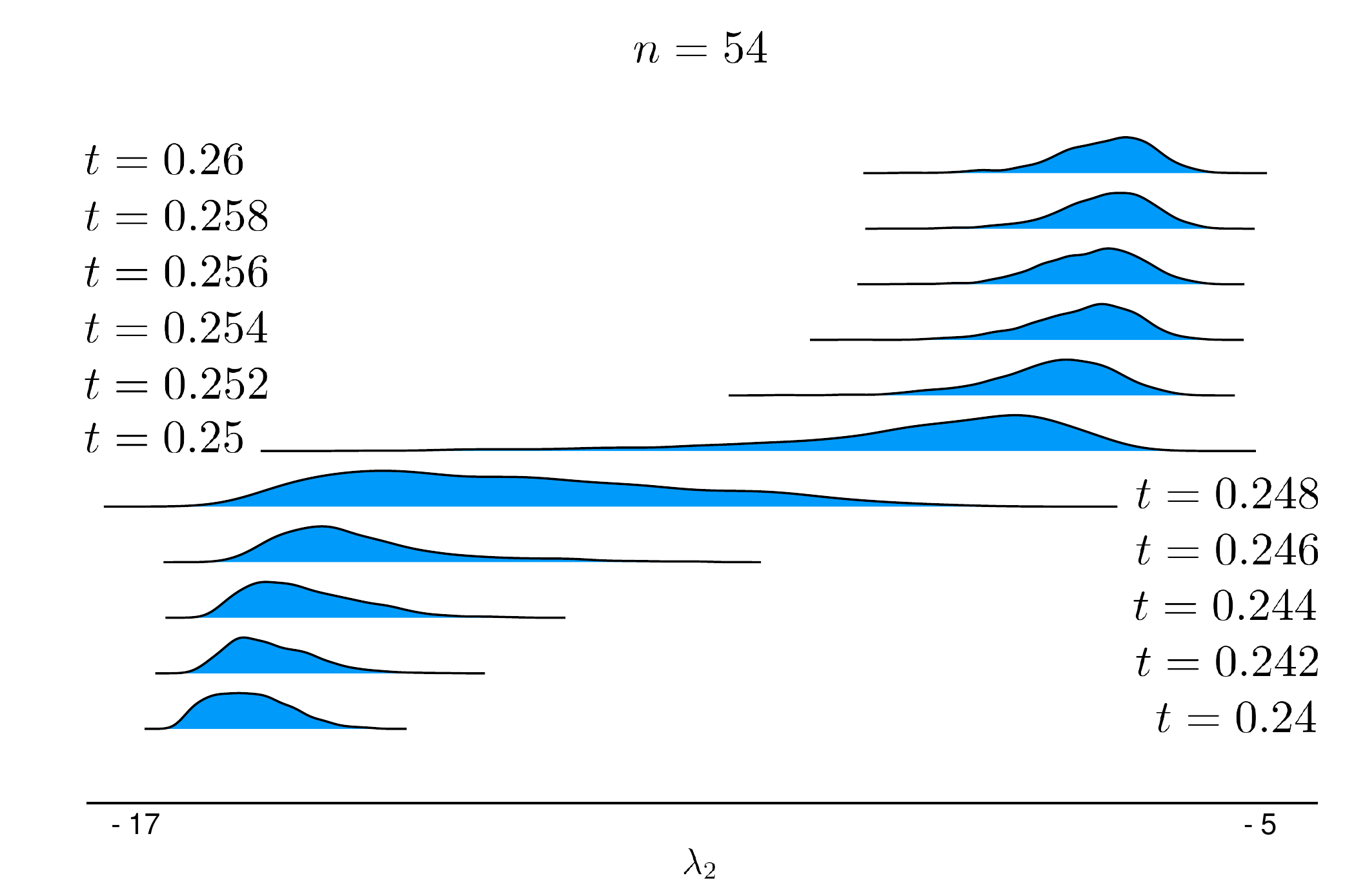}}
    \caption{Density plots showing the distribution of the second most
    negative eigenvalue $\lambda_2$, for graphs with $54$ nodes, edge density
    $0.67$, and various triangle densities $t$.}
\label{histo-evs}
\end{figure}

The histograms for $n=54$, $e=0.67$ and equal jumps of $t$ from $0.24$
to $0.26$ are shown in Figure~\ref{histo-evs}. When $t=0.248$ or
$t=0.250$ the variance in $\lambda_2$ is large, but when $t>0.250$ or
$t<0.248$ the variance is much smaller.  We interpret this as meaning
that for the higher $t$ the networks represent one phase and for lower
$t$ they represent a different phase, and where the variance is large
we have some sort of transition. Note that $\lambda_2$ is much more
negative for values of $t$ below the transition than above, which is
what we would expect from a $B(1,1)/ A(3,0)$ transition. The behavior
of the eigenfunctions is harder to quantify but easy to spot. As noted
earlier, Figures \ref{xiplot-a30} and \ref{xiplot-b11} show the
distributions of $(\xi_1(v), \xi_2(v))$ for two representative sample
networks, one substantially below and one substantially above the
transition.

\begin{figure}[htbp]
\center{\includegraphics[width=5in]{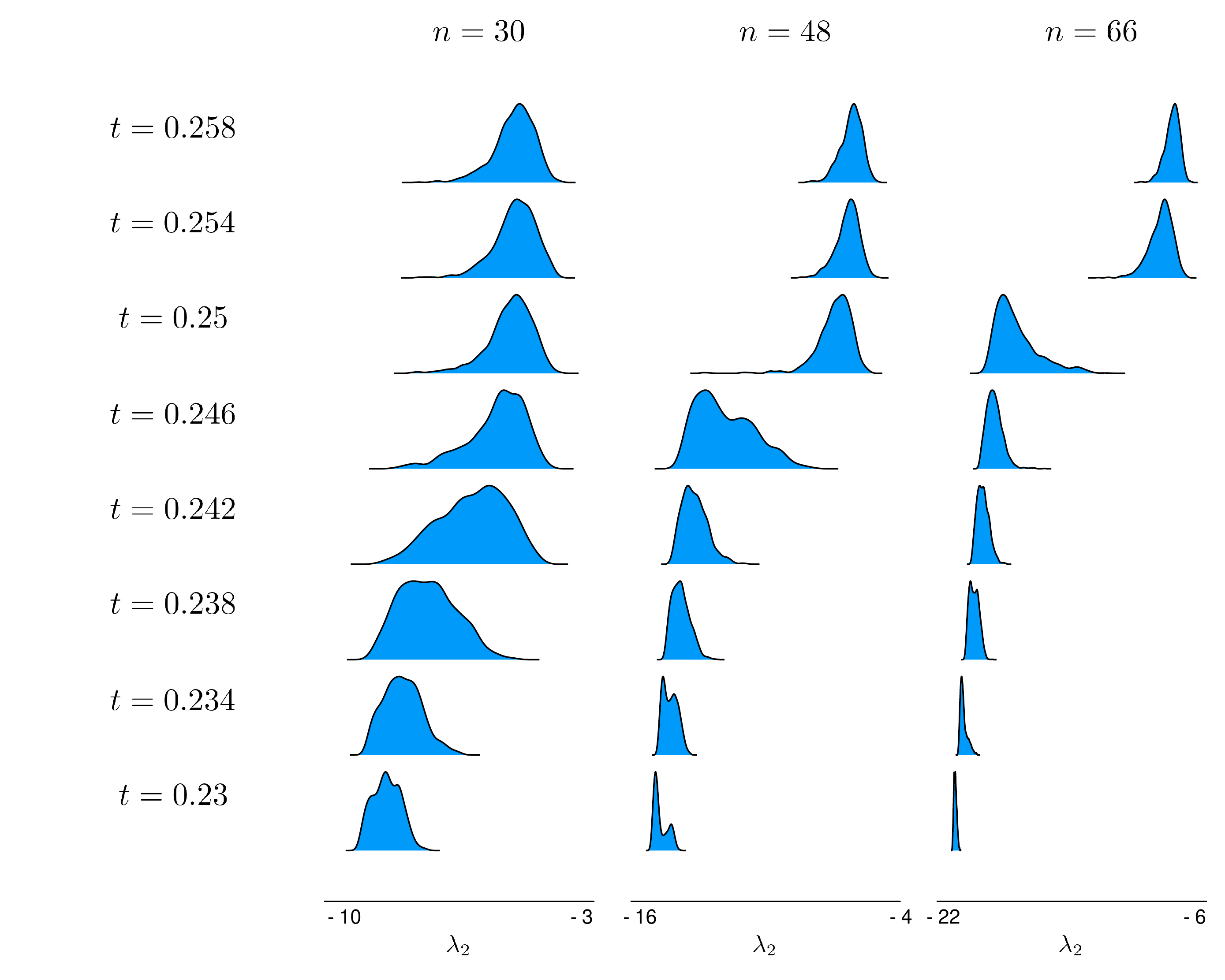}}
\caption{Density plots showing the distribution of the second most negative
    eigenvalue $\lambda_2$ for $n=30$, 48 and 66 nodes, edge density $0.67$,
    and various triangle densities $t$.}
\label{histo-evs-vary-n}
\end{figure}

Figure \ref{histo-evs-vary-n} shows what happens when $n=30$, $48$ or
$66$. Qualitatively, the picture is the same as for $n=54$.
Quantitatively, there are some small differences.  The transition is
sharper when $n$ is larger, but is still visible when $n$ is as small
as 30.  The specific location of the transition varies somewhat in
$n$. As $n$ increases, the transition point approaches from below the
known location of the $B(1,1)/A(3,0)$ phase transition for infinite
$n$, namely $t=0.2596$.

The dependence of the transition $t$ on $n$ is particularly
intriguing.  There is a lower limit to the value of $t$ achievable by
a bipodal network in which each pode has an Erd\H{o}s-R\'{e}nyi structure. 
For the range of values of $n$ considered in this paper
(30 to 66), the transition actually occurs {\em below}
this point! The networks just above the transition have a clear division
into ``red'' and ``blue'' nodes, but the number of red-red-red
triangles is significantly less than one would expect from an
Erd\H{o}s-R\'{e}nyi graph with the given number of red-red edges. 
That is, in this range of $t$-values the red pode
has some sort of $n$-dependent internal structure that we do not
yet understand.  This phenomenon is the subject of continuing
research.

\begin{figure}[htbp]
\center{\includegraphics[width=5in]{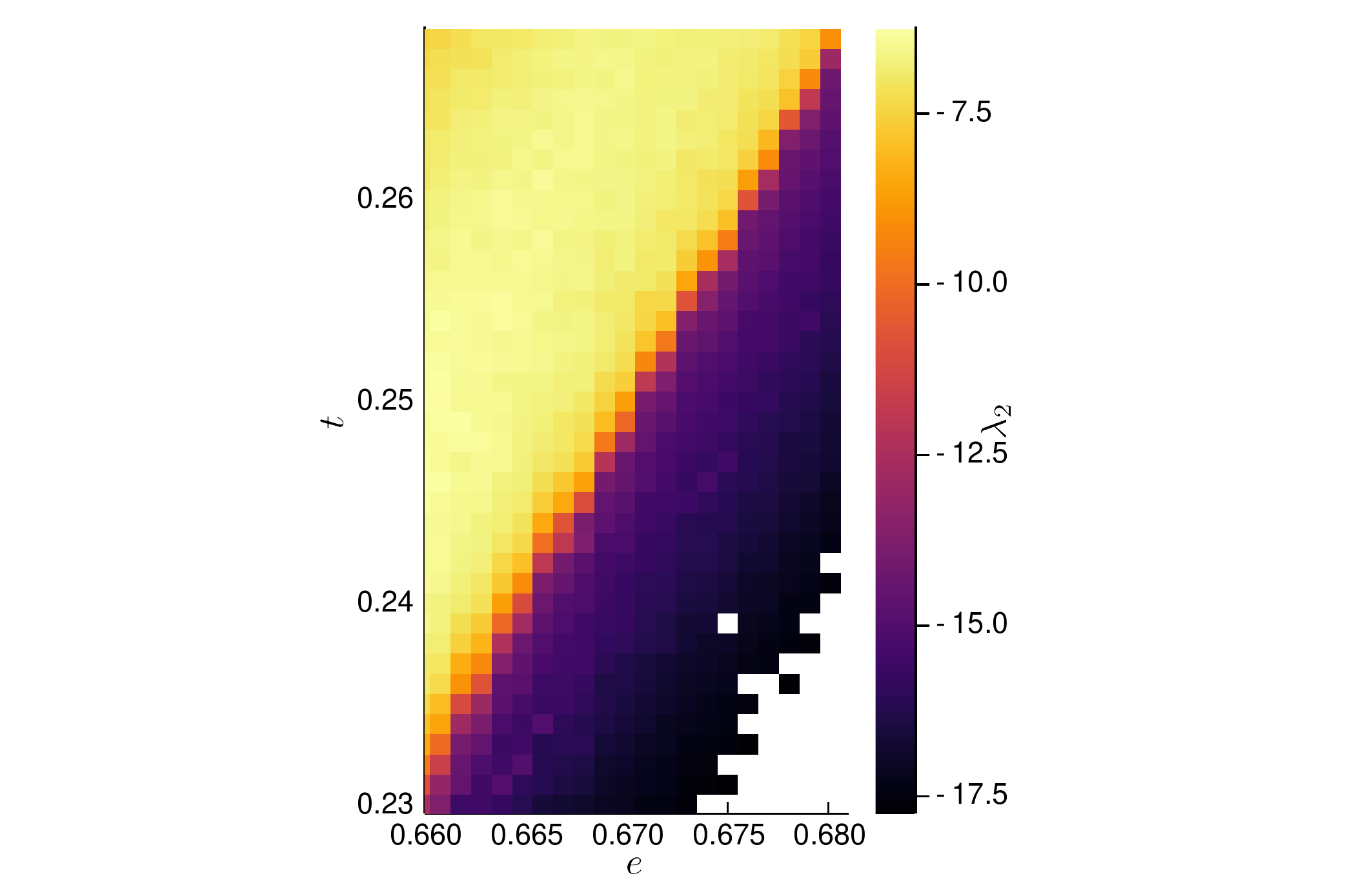}}
\caption{The second-most-negative eigenvalue $\lambda_2$ of a typical
    $54$-node graph, as a function of edge density $e$ and triangle
  density $t$.
    The phase transition is clearly visible as a curve along which
    $\lambda_2$ changes rapidly as a function of $e$ and $t$.}
\label{e-t-plot}
\end{figure}

Finally, we examine whether our results are sensitive to changes in
the edge density $e=0.67$.  They are not. Figure \ref{e-t-plot} shows
the mean value of $\lambda_2$ as a function of $(e,t)$.  There is a
sharp transition for all values of $e$. All that changes with $e$ is the
specific value of $t$ at which it occurs.

Having established that the phenomena we are investigating are
insensitive to changes in $e$ or $n$, we will henceforth restrict
attention to $n=54$, $e=0.67$.

\subsection{Sampling from equilibrium}
\label{sampling}

We use Markov Chain Monte Carlo (MCMC) to (approximately) uniformly sample
graphs with specific numbers of edges and triangles. In order to ensure that
our results are not sensitive to the sampling procedure, we carry this out in
several different ways, all of which follow the same
general procedure: given fixed number $n$ of vertices, and given fixed edge and
triangle counts $E$ and $T$, we define the following three objects:
\begin{enumerate}
    \item a constraint set $\Omega$ which contains all graphs on $n$
        vertices with $E$ edges and $T$ triangles,
    \item a function $d_\Omega$ measuring the ``distance'' from a graph to
        $\Omega$, and
    \item for each graph $G$, a probability measure $Q_G$ on graphs with $n$
        nodes (the ``proposal distribution'').
\end{enumerate}

We use this data to define a Metropolis-Hastings Markov chain with the
following update rule: given the current state $G$ (a graph on $n$ nodes), we
draw the graph $G'$ from $Q_G$. If $G' \in \Omega$, it becomes the next state
of the Markov chain. If $G' \not \in \Omega$ but $d_\Omega(G') \le
d_\Omega(G)$, $G'$ becomes the next state of the Markov chain. Otherwise, the
Markov chain remains at $G$.

The Markov chain defined above is reversible, and the uniform measure on
$\Omega$ is a stationary measure (and if $\Omega$ is connected under the
proposal distributions, this stationary measure is unique). Finally, we
generate a random graph by running the Markov chain for a long time (for most
of the empirical results presented here, $10^9$ accepted steps were sufficient)
and accepting the final graph if it had exactly $E$ edges and $T$ triangles.

All of the results presented here were checked with the
following three choices of $\Omega$, $d_\Omega$, and $Q$:
\begin{enumerate}
    \item $\Omega$ is the set of graphs with exactly $E$ edges and at most $T$
        triangles; $d_\Omega(G) = \max\{0, T(G) - T\}$; and $Q_G$ is the random
        graph obtained from $G$ by deleting an edge at random and adding an
        edge at random
    \item $\Omega$ and $d_\Omega$ are as above;
        $Q_G$ is the random graph obtained from $G$ by choosing a vertex at random,
        then deleting an incident edge at random and adding an incident edge at random
    \item $\Omega$ is the set of graphs with $|E(G) - E| \le C$ and $|T(G) - T| \le C$
        for some constant $C$; $d_\Omega(G) = \max\{0, |E(G) - E| + |T(G) - T| - 2C\}$;
        and $Q_G$ is the random graph obtained from $G$ by either deleting a random
        edge, adding a random edge, or deleting a random edge and adding a random edge.
\end{enumerate}

The figures and numbers that we present in the main body of this work are all
for the first choice of dynamics. The other two choices of dynamics showed
the same qualitative behavior but mixed more slowly; the third choice of
dynamics in particular has a tendency to ``get stuck'' for long periods, since
after deleting an edge it sometimes struggles to find another move that
respects the triangle constraints.

\subsubsection{Evaluating convergence}\label{subsubsection}

Since we are not able to prove bounds on the mixing time of any of the Markov
chains above, we validated the convergence of our Markov chains by verifying
that the distribution of certain statistics of interest (for example, the
second-smallest eigenvalue of the resulting graphs) were independent of the
Markov chains' initial states. Specifically, we initialized our various Markov
chains from three different initial distributions:
\begin{enumerate}
    \item a uniformly random graph with exactly $E$ edges,
    \item a graph drawn from a $B(1,1)$ graphon, and
    \item a graph drawn from an $A(3,0)$ graphon.
\end{enumerate}
Note that these three distributions have very different spectral properties.
After running the Markov chains above, however, we verified that the spectral
statistics of the outputs did not depend on which of the three initializations
we chose.

\section{Structural rearrangement under change 
of phase}\label{nucleation}

We now turn to the central question of this paper. We have seen that
for $t=t_{low}=0.24$, the ensemble is dominated by $A(3,0)$ networks,
as evidenced by $\lambda_2$ being very negative, and by
$(\xi_1,\xi_2)$ being clustered around three points, with equal
numbers of nodes near each point. When $t=t_{high}=0.26$, the
ensemble is dominated by networks with $\lambda_2$ much smaller, with
the $(\xi_1,\xi_2)$ plot being a vertical bar (representing 60\% of
the nodes) and a separate cluster (40\%). But how does the network
change from one kind of network to the other under our dynamics, 
described in Section~\ref{equilibrium}?

The supplementary material contains several
\href{https://www.youtube.com/playlist?list=PLZcI2rZdDGQqx3WEY6BXoJXcxqQXF5ZzQ}{animations}
that show the evolution
of $(\xi_1, \xi_2)$ and various other statistics across the three stages. In
this section, we do our best to summarize the main features using figures and
text. However, we highly recommend viewing the actual animations.
A picture is worth a thousand words, and a movie is worth a thousand 
pictures. 

\subsection{Increasing $t$}

\begin{figure}[htbp]
\center{\includegraphics[width=5in]{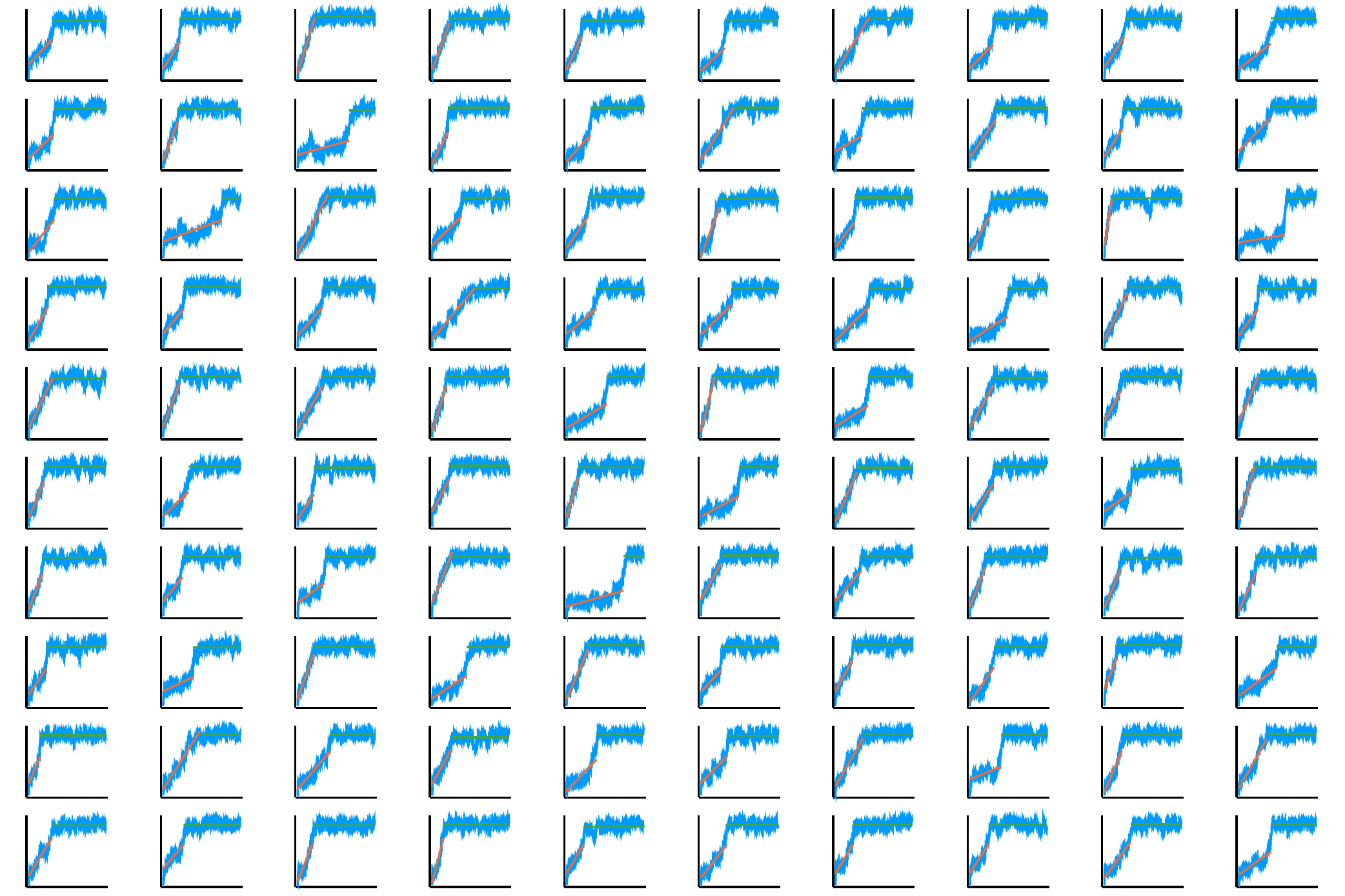}}
\caption{100 repetitions showing the trajectory of $\lambda_2$ as a function of time, as the triangle
    density is increased from 0.24 to 0.26 over 10,000 steps. The red and green lines show the least-squares best fit
    of the trajectory by a linear function (the red line) followed by a constant function (the green line).
    }
\label{short-up}
\end{figure}

\begin{figure}[htbp]
\center{\includegraphics[width=4in]{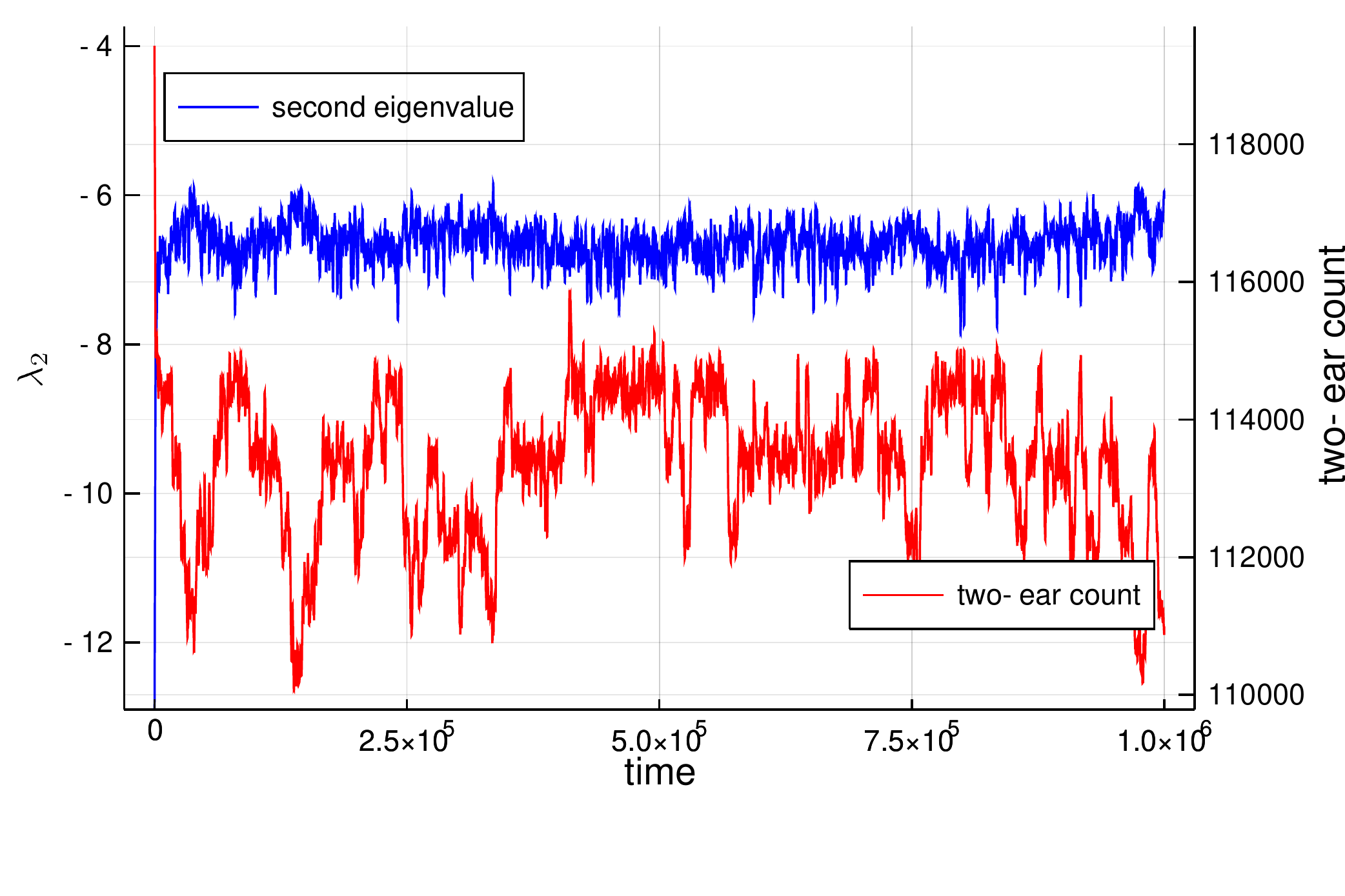}}
\center{\includegraphics[width=4in]{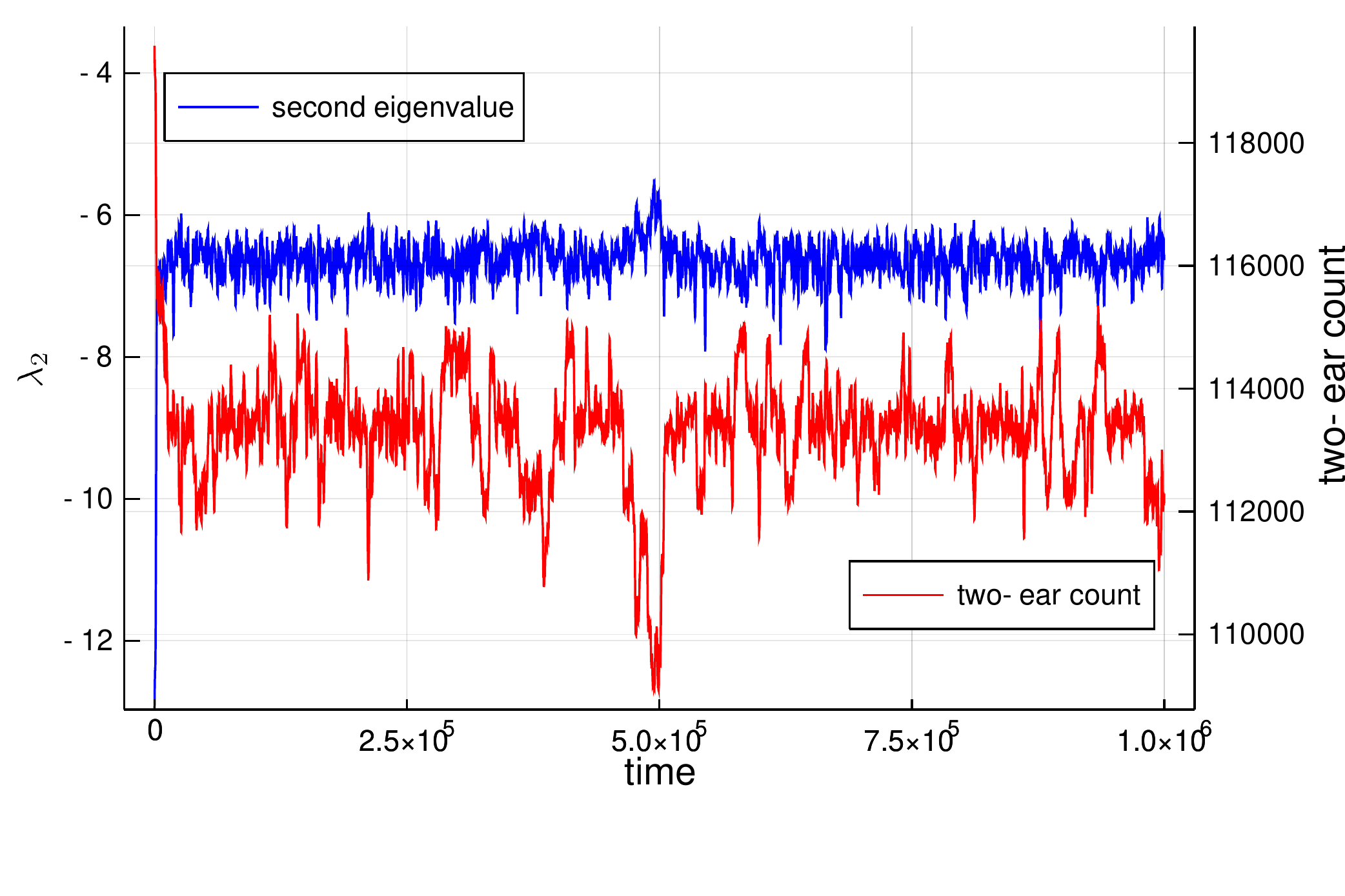}}
    \caption{Two different random trajectories of $\lambda_2$ (in blue), and the number of
    embedded ``two ear'' graphs (in red) as a function of time, as the triangle
    count is increased from 0.24 to 0.26 over $10^6$ steps.
    The data is smoothed by averaging over windows of 1000 steps.
    }
\label{long-up}
\end{figure}

We display 100 short runs of 10,000 steps in Figure~\ref{short-up}
and two long runs of 1,000,000 steps in Figure~\ref{long-up}.
Both figures
show $\lambda_2$ as a function of 
time. Figure~\ref{long-up} also shows the number of embedded
``two ear'' graphs, where a ``two ear'' is two triangles sharing a 
common edge, or equivalently a tetrahedron with one edge missing. 
The data in Figure~\ref{long-up} has been smoothed, with both 
$\lambda_2$ and the 2-ear count averaged over 1000 successive times.
The data in Figure~\ref{short-up} has not been smoothed. 

A typical transition occurs in three dynamical stages. The first two
are visible in Figure~\ref{short-up} and the third is visible in
Figure~\ref{long-up}.  In the first stage, which only takes 
about 50 steps, the network retains its $A(3,0)$ structure as $t$
increases from $t_{low}$ to $t_{high}$. In this stage $\lambda_2$
increases somewhat as the parameters $a$ and $b$ change, but the
distribution of $(\xi_1, \xi_2)$ does not significantly change.
This stage is visible in Figure~\ref{short-up} as a sudden, but modest,
increase in $\lambda_2$ at the very beginning of each plot.

In the second stage, which usually takes between 1000-6000 steps, the
network loses its $A(3,0)$ structure as two of the clusters
merge. This stage is marked by an increase in $\lambda_2$ and a
dramatic change in the distribution of $(\xi_1, \xi_2)$. At the end of
this stage, we have a structure that is very similar to the optimal
$B(1,1)$ graphon, except that the podes do not have the optimal sizes,
being closer to $2:1$ than to 60:40. This stage is visible in
Figure~\ref{short-up}, and accounts for most of the visible
change in $\lambda_2$. In the last stage, which takes
hundreds of thousands of steps, the relative sizes of the podes adjust
back and forth among a number of possibilities with nearly equal
entropies, resulting in fluctuations in $\lambda_2$ and the 2-ear
count over very long time scales. The beginning of this stage can be
seen in Figure~\ref{short-up} as the point at which $\lambda_2$ becomes
approximately constant over time. In order to actually see the pode
sizes changing, one has to view the trajectories at the larger
time-scale shown in Figure~\ref{long-up}.

The first stage is driven entirely by entropy.  Since all networks are
well below the target number of triangles, all swaps of one edge for
another are allowed. This results in edges simply moving from where
they are concentrated to where they are not. Since $a<b$, $a$ gets
bigger and $b$ gets smaller, but (aside from random noise) the graphon
remains piecewise constant, with the three podes remaining
indistinguishable. This is reflected in the fact that plots of
$(\xi_1, \xi_2)$ basically do not change. Since $\lambda_2$ is
proportional to $a-b$, $\lambda_2$ increases somewhat. For $n=54$, it
goes from a little below $-16$ to around $-14.5$.

\begin{figure}[htbp]
    \includegraphics[width=2in]{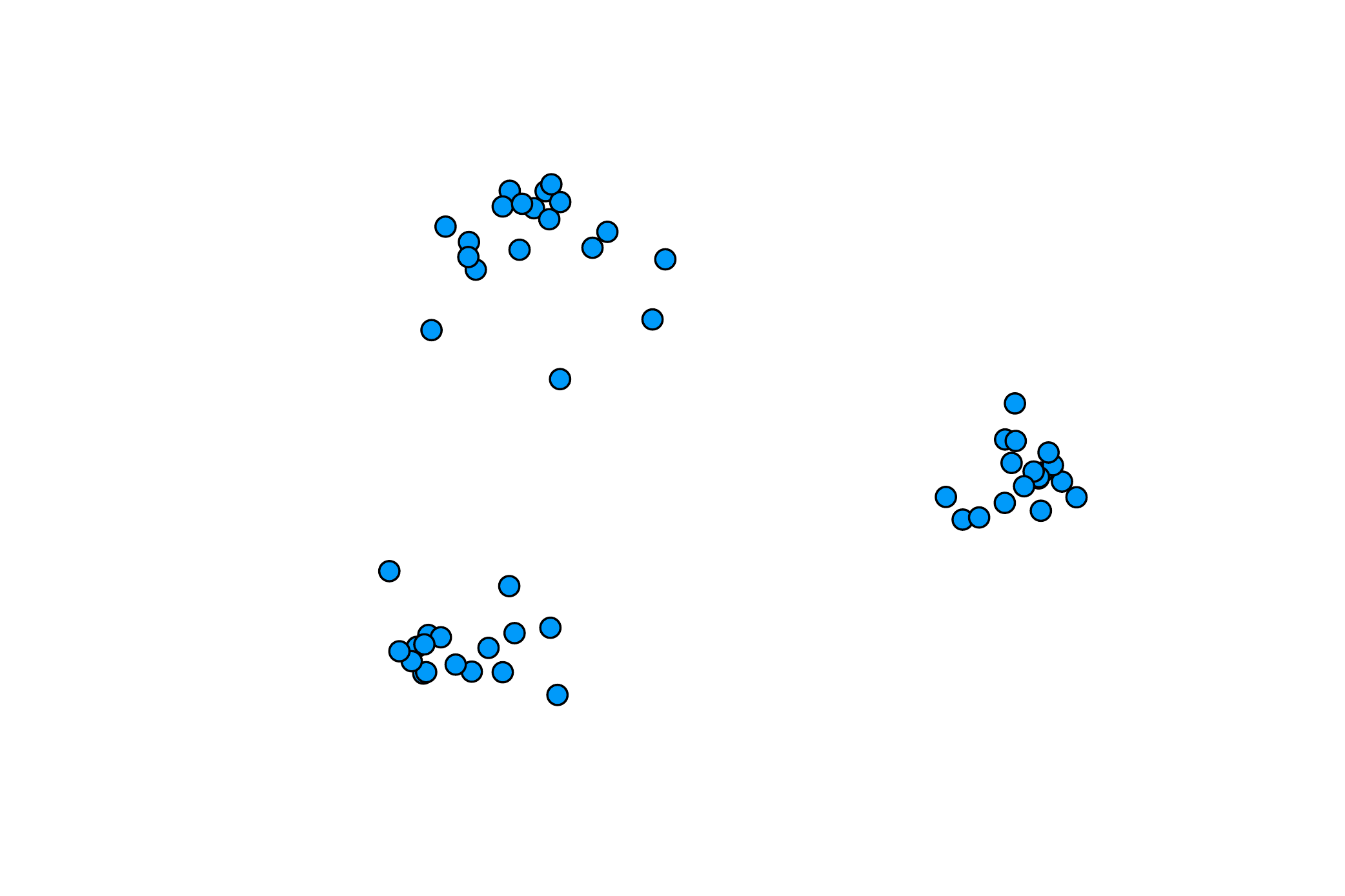}
    \includegraphics[width=2in]{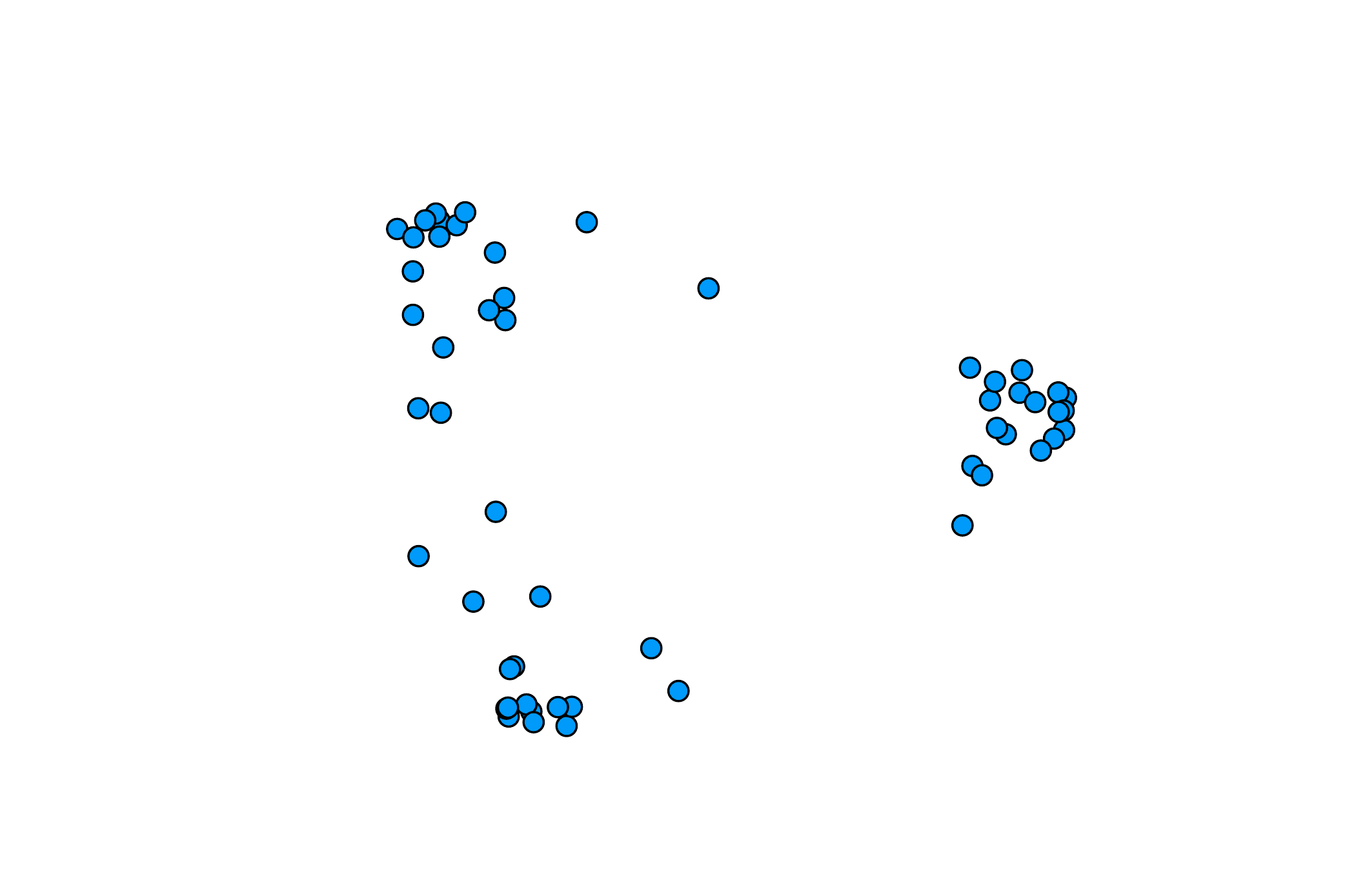}
    \includegraphics[width=2in]{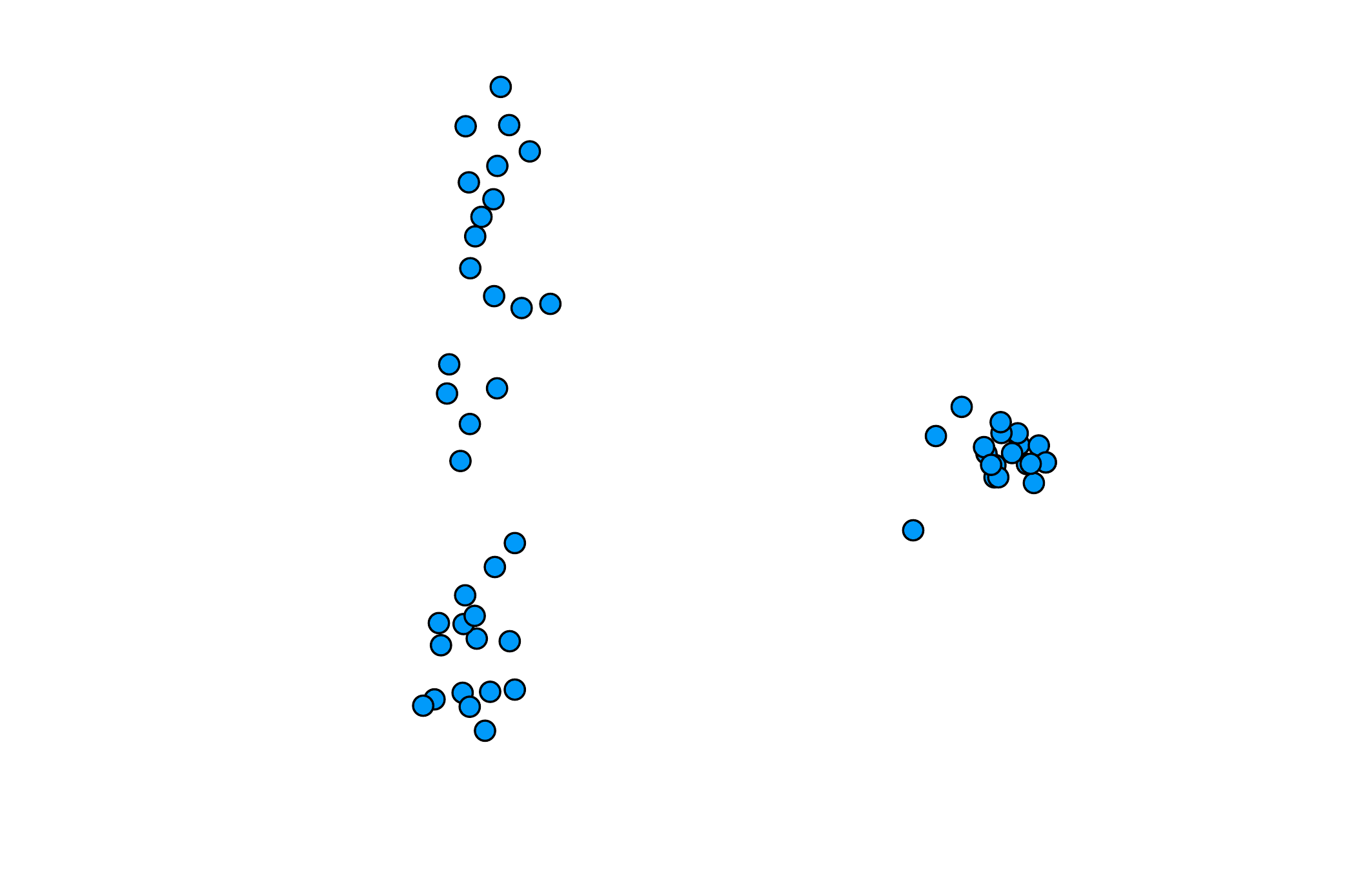}
\caption{Snapshots of $(\xi_1, \xi_2)$ after 500, 800, and 1500 
steps during an upward trajectory. Each dot corresponds to a vertex $v$;
the coordinates of the dot are given by $\xi_1(v)$ and $\xi_2(v)$.
The transition can also be viewed
in the \href{https://www.youtube.com/playlist?list=PLZcI2rZdDGQqx3WEY6BXoJXcxqQXF5ZzQ}{animations}.}
\label{Xi-up}
\end{figure}

In the second stage, the network explores the available phase space
for $t=t_{high}$ starting from an $A(3,0)$ structure. Within each
cluster, the variances in $\xi_1$ and $\xi_2$ increase, and an
occasional vertex may move from one cluster to another.
In a typical run, there is a
latency period 
during which $\lambda_2$ gradually increases from around $-14.5$ to
around $-11$, after which two of the clusters join and $\lambda_2$ rises
quickly to around $-7$. During the merger, all three
clusters spread out, and sometimes the third cluster gains or loses a
vertex or two. Figure~\ref{Xi-up} shows the evolution of $\xi_1$ and $\xi_2$
during this second stage.

The third stage involves punctuated equilibrium. Most of the time,
both $\lambda_2$ and the 2-ear density fluctuate in a range, and the
distribution of $(\xi_1,\xi_2)$ does not change qualitatively. Once in
a while, a node moves from one pode to another, triggering sudden
changes in both $\lambda_2$ and the two-ear count.
When the new ratio of pode sizes is close to the optimal ratio, the system
can stabilize in the new state for almost 100,000 steps. If the new
ratio of pode sizes is far from optimal, however, entropy drives them to
change back more quickly, usually within 10,000 steps. Such an extreme excursion
is visible at time 150,000 in the first run, where the smaller pode grows
to 24 nodes and remains there for about 10,000 steps before eventually
making its way back to 21 nodes at time 170,000. On its way back, the
size of the smaller pode briefly (for about 10,000 steps) stabilizes
at 23, corresponding to the two-ear count stabilizing at around 111,150.

The end of stage 1 and the beginning of stage 2 is easy to
define. This is when $t$ first hits its target value, after which the
constraint on $t$ begins to affect which edge swaps are allowed. This
consistently takes between $45$ and $60$ steps, with an average of $50$.

\begin{figure}[htbp]
\center{\includegraphics[width=4in]{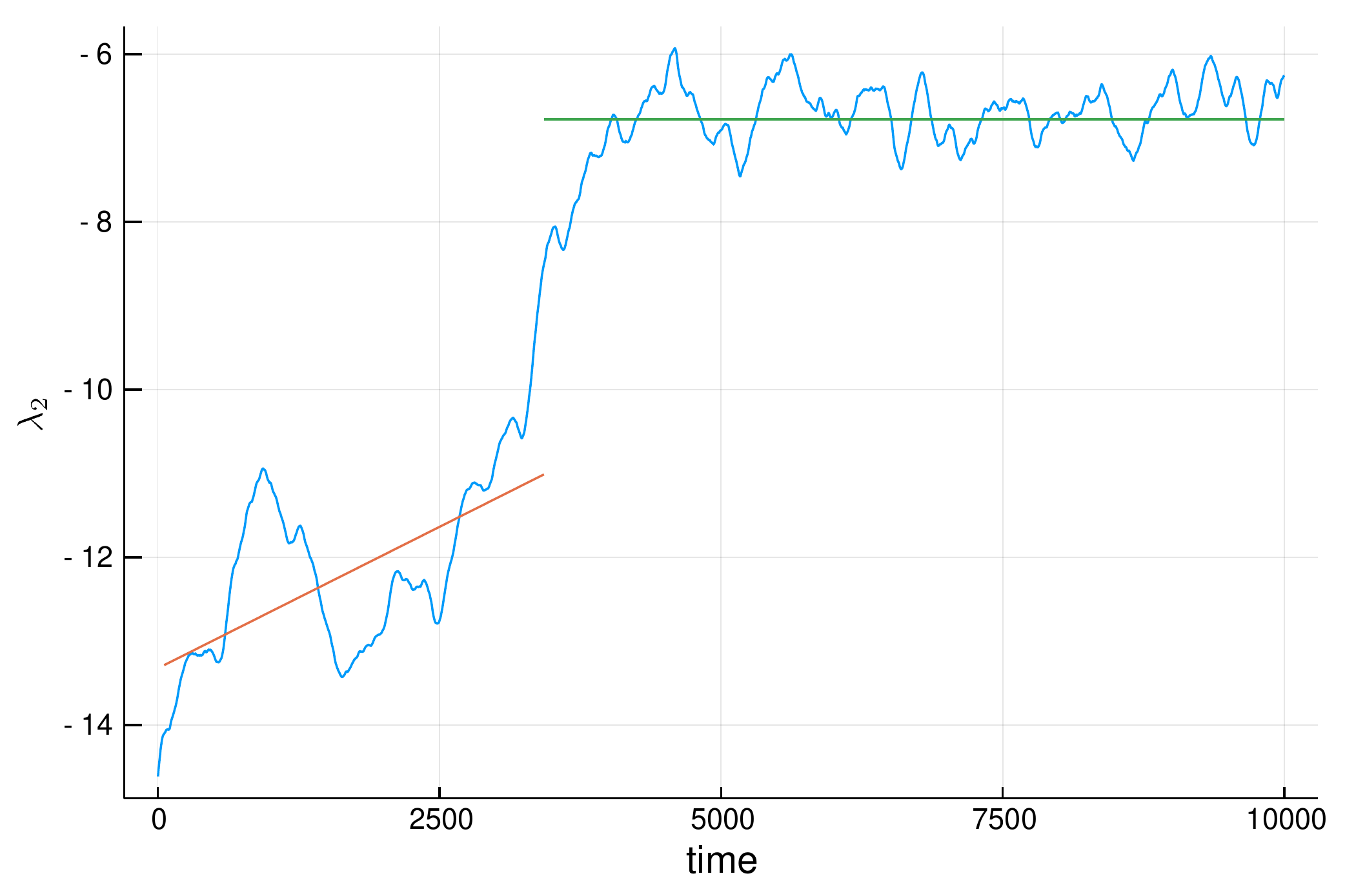}}
\center{\includegraphics[width=4in]{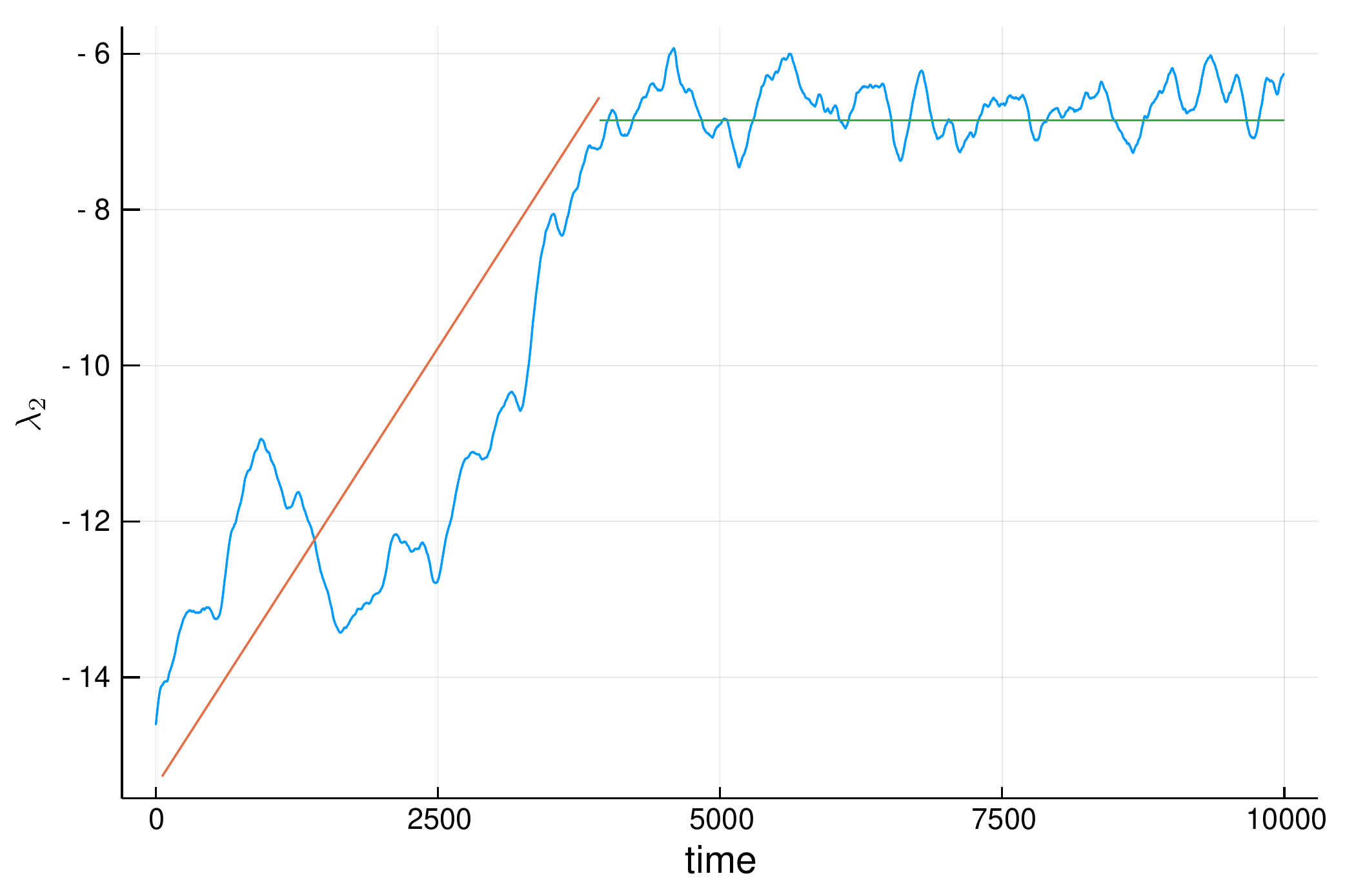}}
\center{\includegraphics[width=4in]{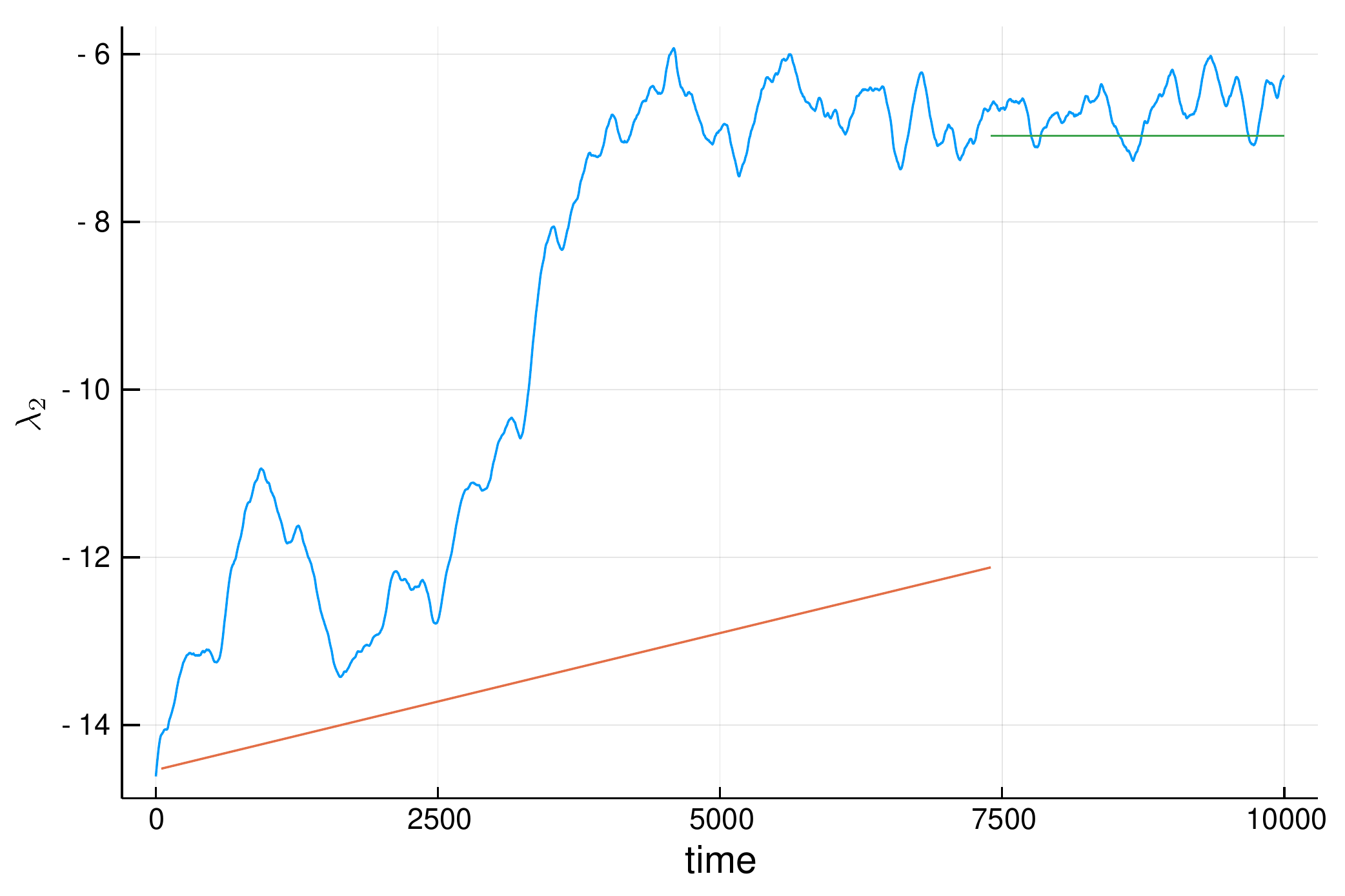}}
    \caption{Close-ups of the 17th, 23rd, and 65th trajectories from Figure~\ref{short-up},
    smoothed over 100 steps.
These three trajectories can also be seen in the
\href{https://www.youtube.com/playlist?list=PLZcI2rZdDGQqx3WEY6BXoJXcxqQXF5ZzQ}{animations}.}
\label{3-short-up}
\end{figure}

The end of stage 2 / beginning of stage 3 is harder to define. To do
this, we did a least-squares fit of the first 10,000 steps, minus
stage 1, via a line of arbitrary slope followed by a horizontal
line. These best fit lines are shown in Figure~\ref{short-up}.  We
interpret the domain of the sloped line as stage 2 and the domain of
the horizontal line as stage 3.

The evolution of $\lambda_2$ (smoothed over 100 steps) in runs 17, 23
and 65 is shown in more detail in Figure~\ref{3-short-up}.  In a run
with a long latency, such as run 65, the sloped line tracks the
gradual increase during the latency period, and the subsequent rapid
rise shows up as a discontinuity between the lines.  In runs with
short latencies, such as run 23, the entire increase from $-14.5$ to
$-7$ is captured by the sloped line, so there is no discontinuity. In
a few runs, such as run 17, the latency period involves some false
starts, where $\lambda_2$ first increases and then decreases, and the
fit to a sloped line is not good at all. However, even in these
exceptional cases, the beginning of the horizontal line provides a
credible measure of the beginning of stage 3.

\begin{figure}[htbp]
\center{\includegraphics[width=5in]{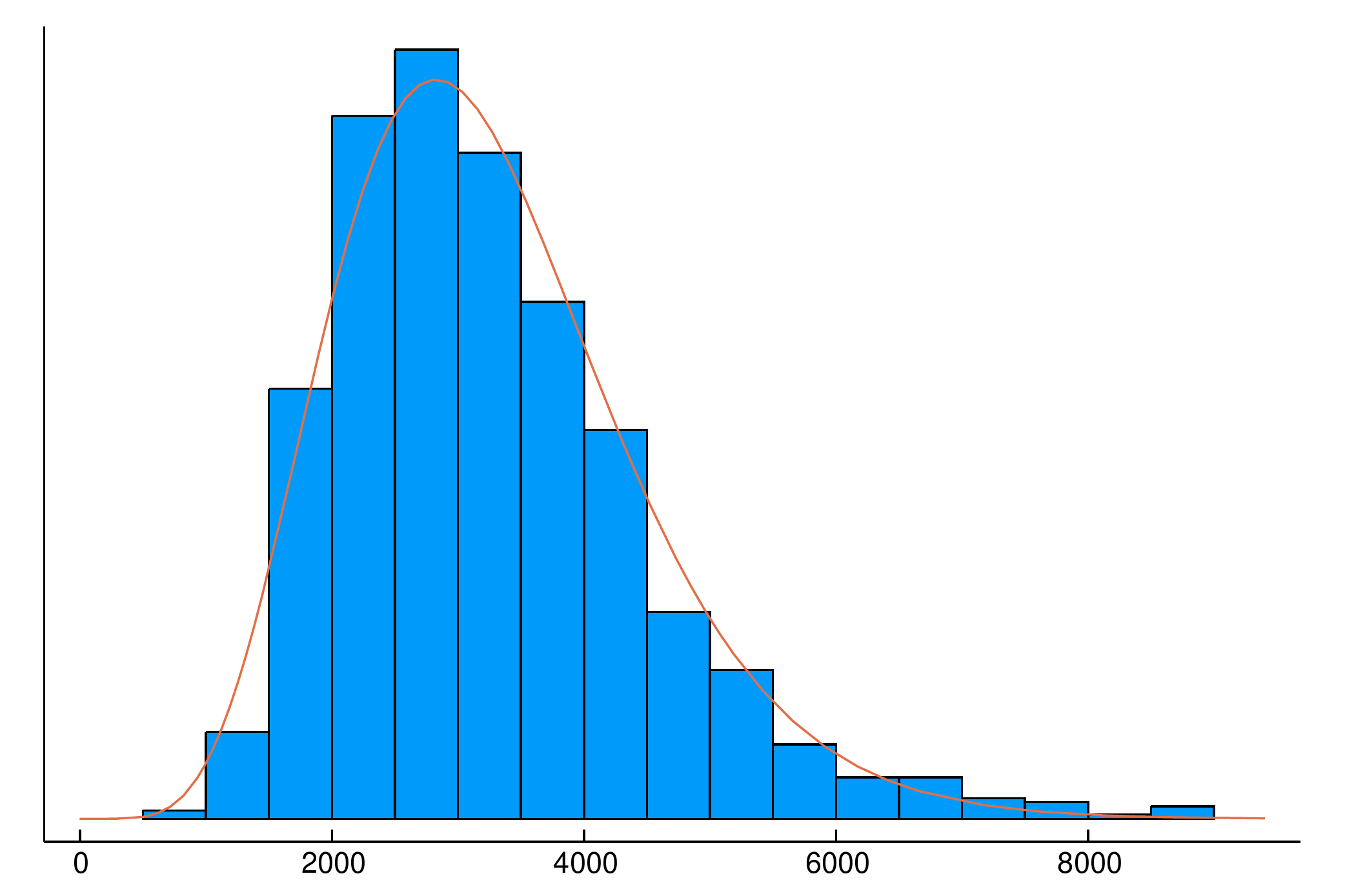}}
\caption{The histogram of the length of stage 2 for upwards trajectories, and a
best-fit Gamma function with ($\alpha = 7.5$ and $\theta = 430$)}
\label{Gamma-up}
\end{figure}

We have plotted the distribution of the length of stage 2 in
Figure~\ref{Gamma-up}. The shape of this histogram is modeled well by
a Gamma distribution with $\alpha=7.5$ and $\theta=430$
(i.e. the
function $x^{6.5} e^{-x/430}$), which we have superimposed on the
histogram.

\subsection{Decreasing $t$}

\begin{figure}[htbp]
\center{\includegraphics[width=5in]{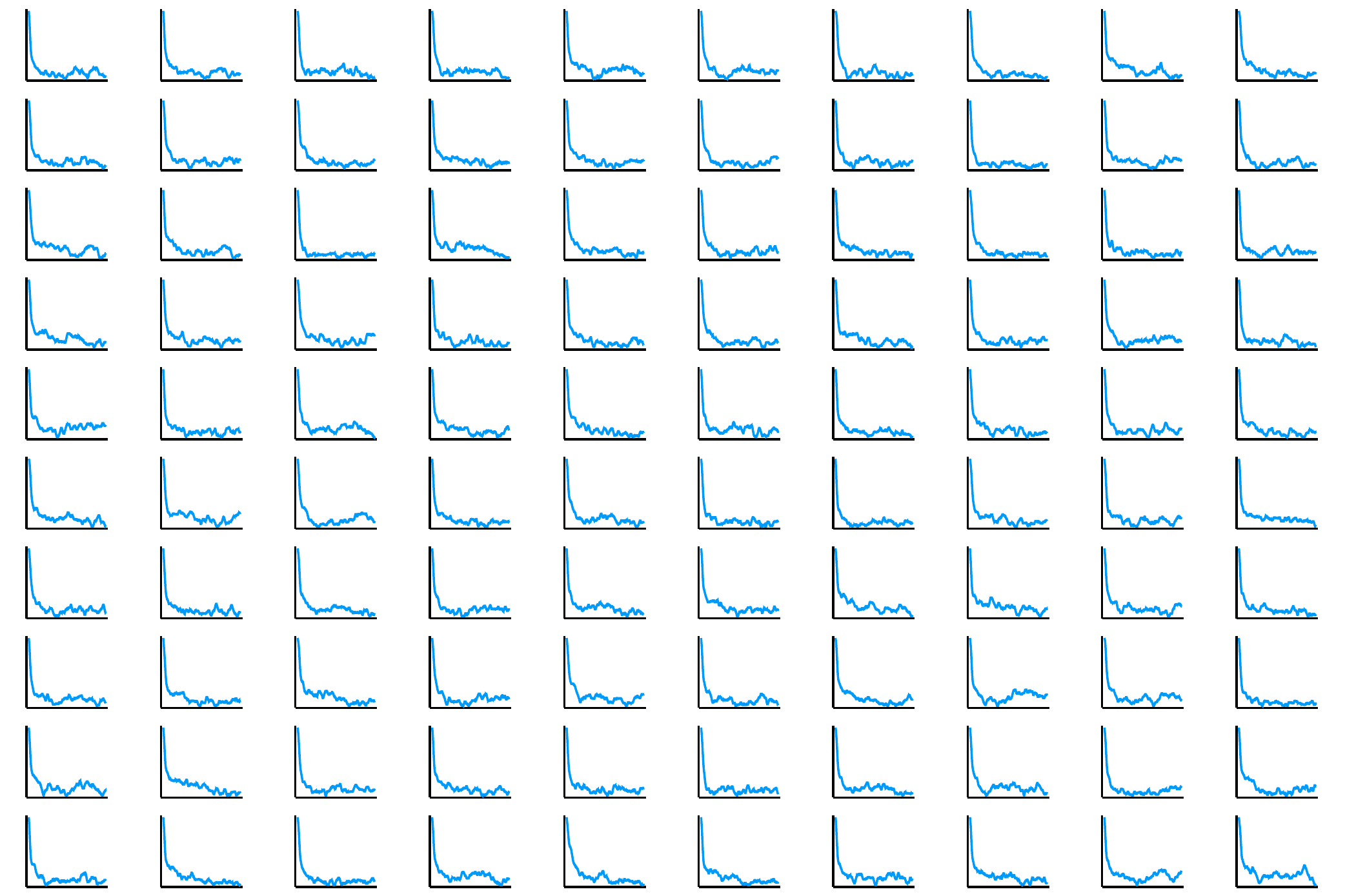}}
\caption{100 repetitions showing the trajectory of $\lambda_2$ as a function of time, as the triangle
    density is decreased from 0.26 to 0.24 over 5000 steps. The red and green lines show the least-squares best fit
    of the trajectory by a linear function (the red line) followed by a constant function (the green line).
    The trajectories are smoothed over 100 timesteps.
    }
\label{short-down}
\end{figure}

\begin{figure}[htbp]
\center{\includegraphics[width=4in]{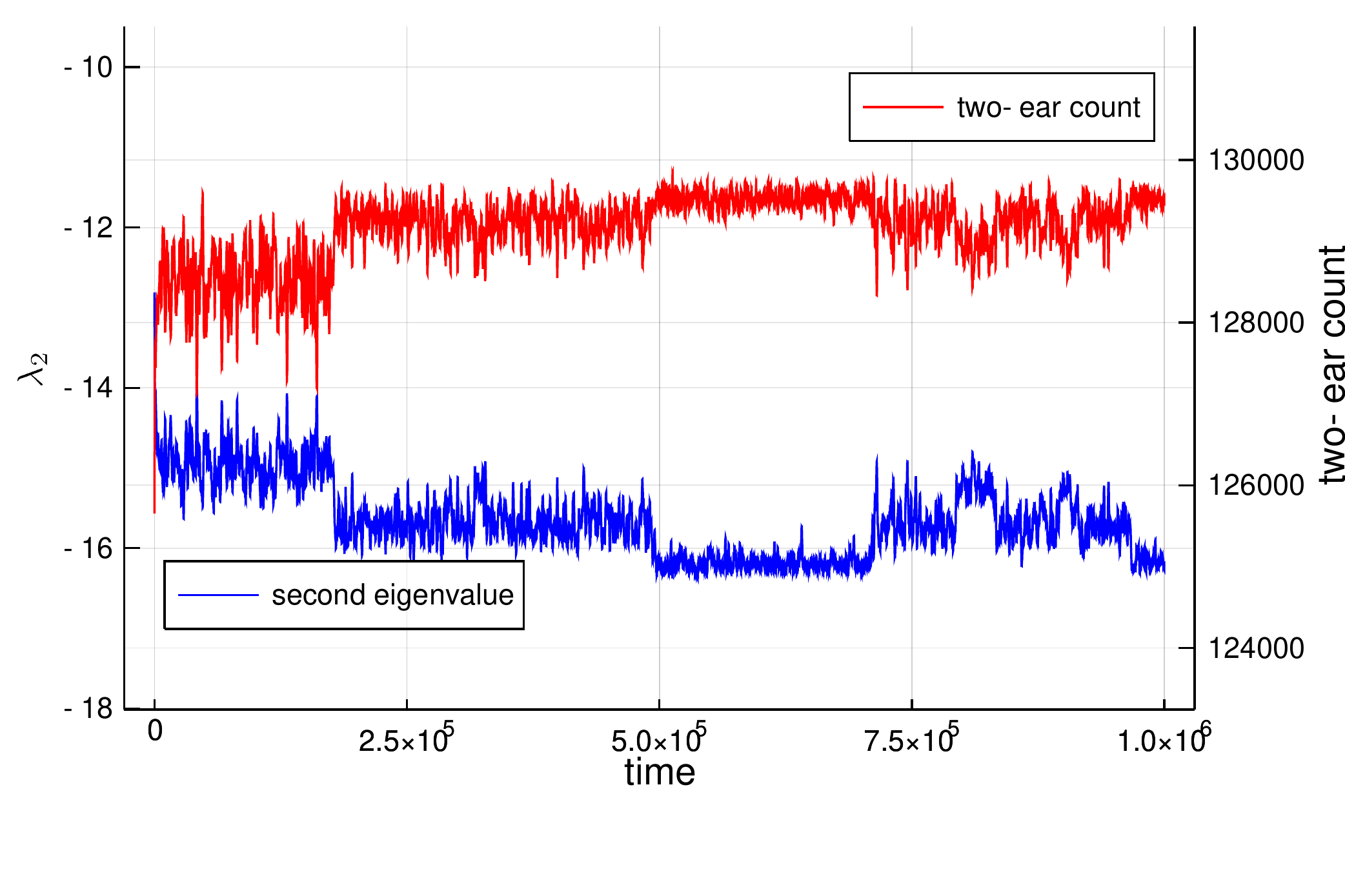}}
\center{\includegraphics[width=4in]{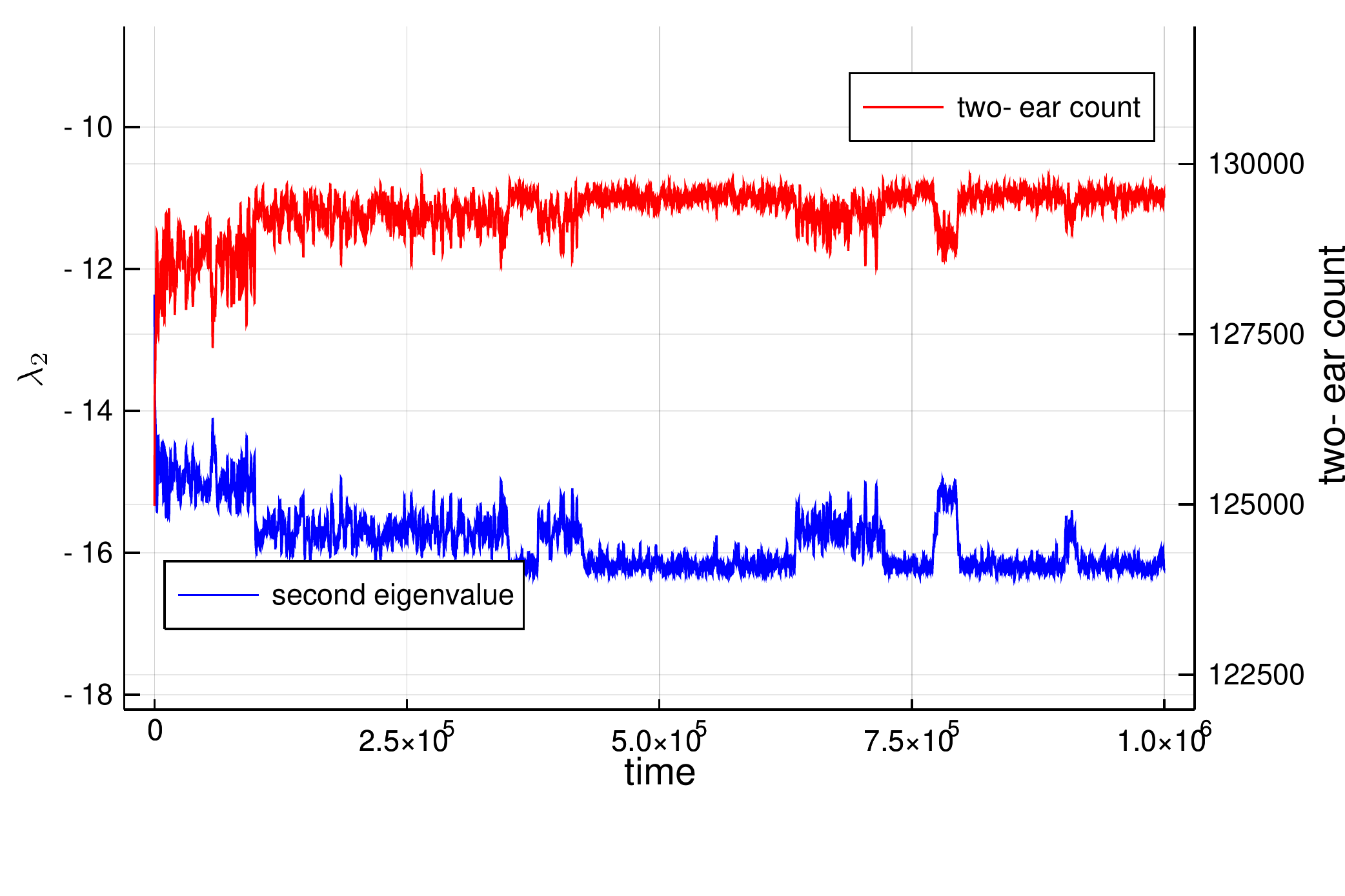}}
    \caption{Two different random trajectories of $\lambda_2$ (in blue), and the number of
    embedded ``two ear'' graphs (in red) as a function of time, as the triangle
    count is decreased from 0.26 to 0.24 over $10^6$ steps.
    The data is smoothed by averaging over windows of 1000 steps.
    }
\label{long-down}
\end{figure}

\begin{figure}[htbp]
\center{\includegraphics[width=4in]{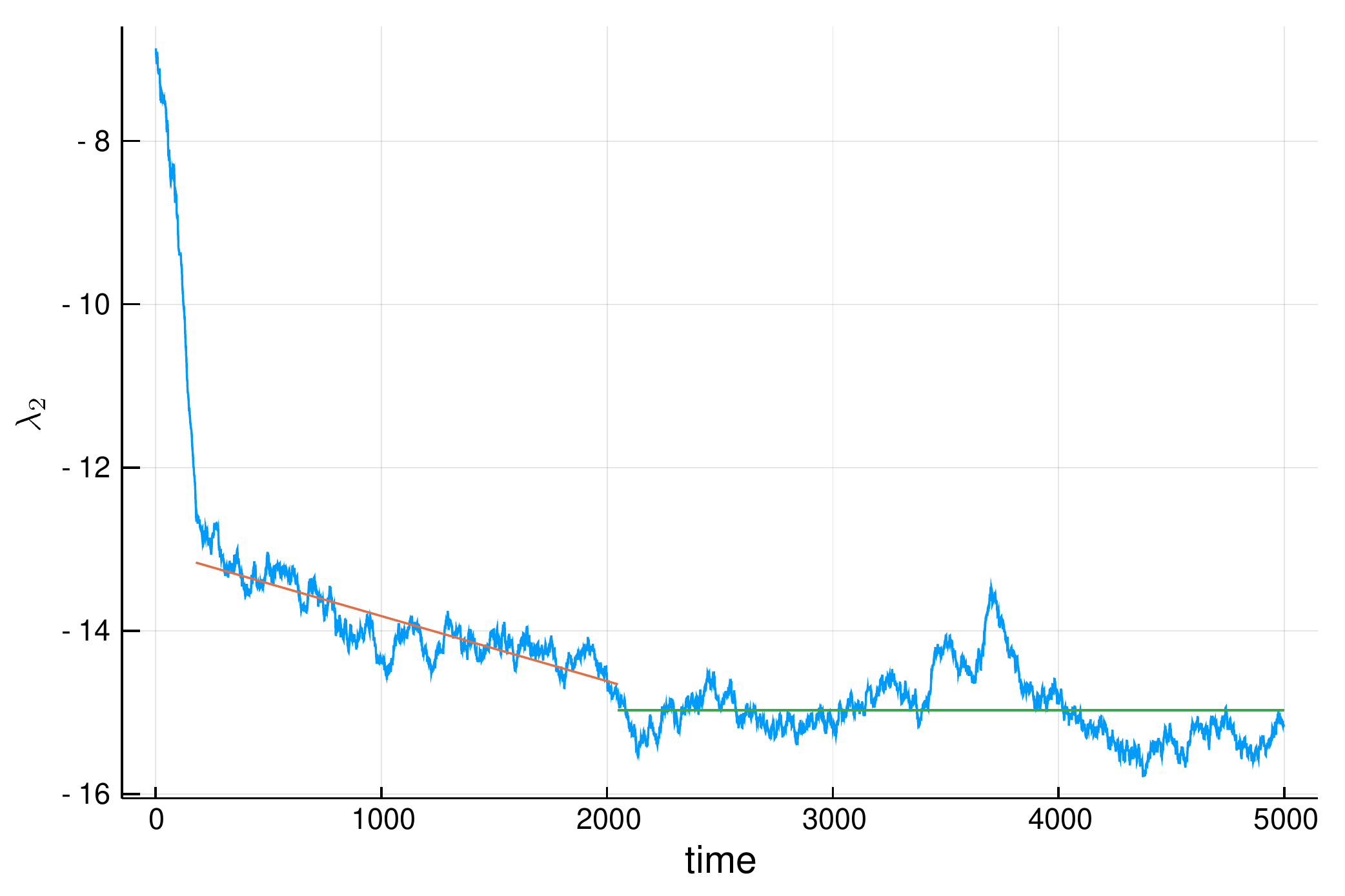}}
\center{\includegraphics[width=4in]{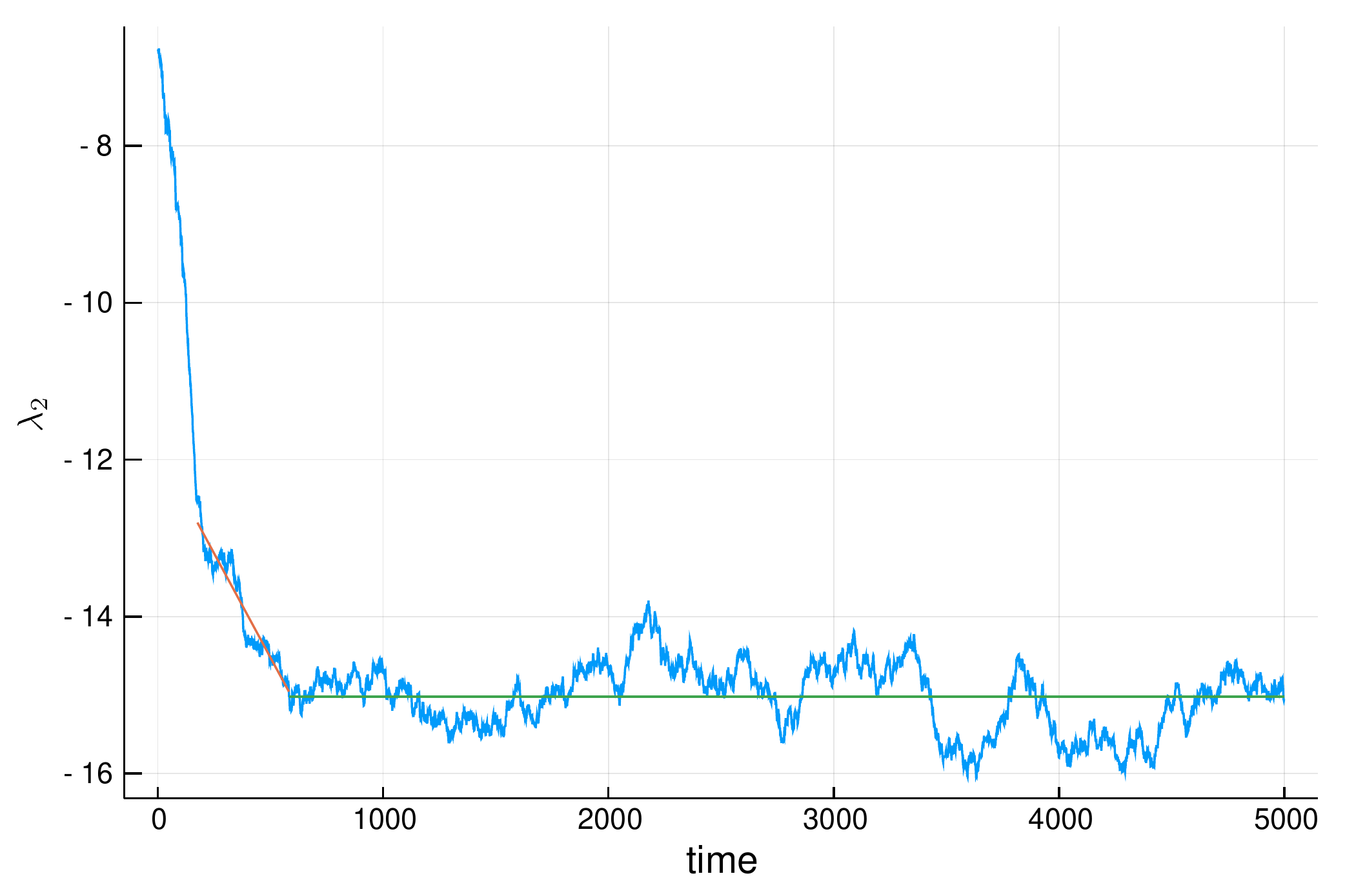}}
\center{\includegraphics[width=4in]{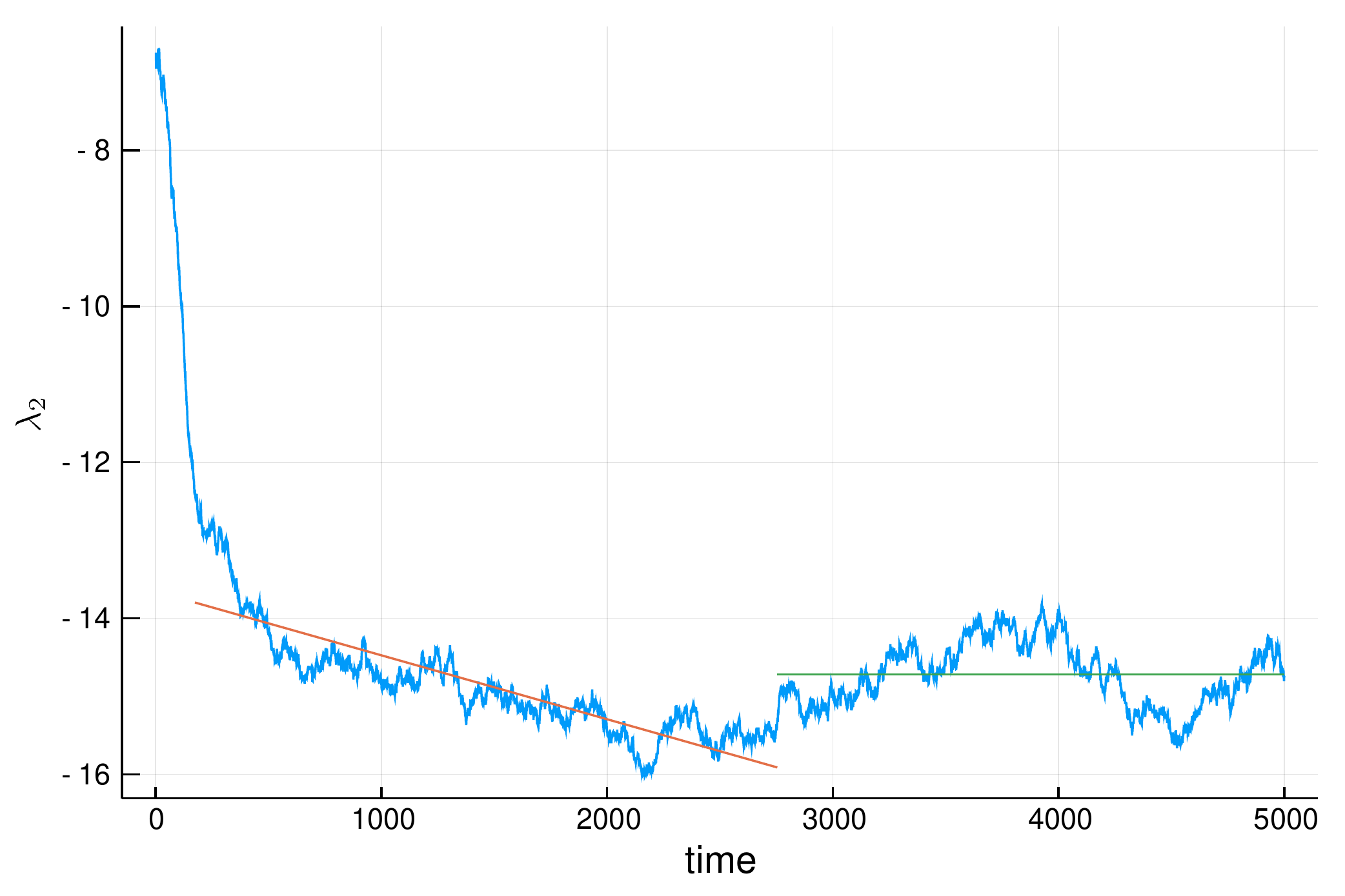}}
    \caption{Close-ups of the 9th, 46th, and 99th trajectories in Figure~\ref{short-down},
    not smoothed.
These three trajectories can also be viewed in the
\href{https://www.youtube.com/playlist?list=PLZcI2rZdDGQqx3WEY6BXoJXcxqQXF5ZzQ}{animations}.}
\label{3-short-down}
\end{figure}

\begin{figure}[htbp]
\center{\includegraphics[width=5in]{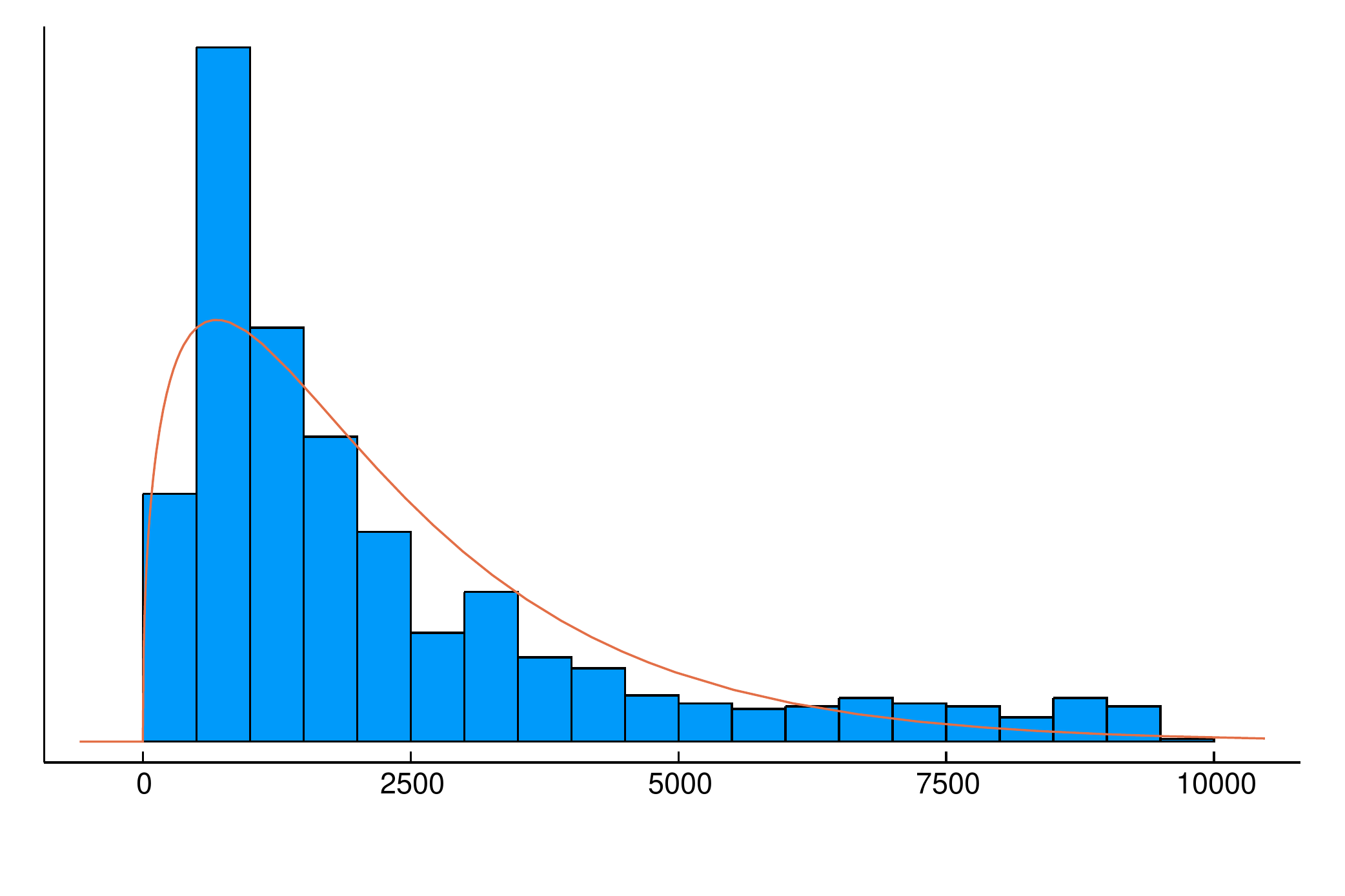}}
\caption{The histogram of the length of stage 2 for downwards trajectories, 
and a best-fit Gamma function with parameters $\alpha=2.25$ and $\theta=660$.}
\label{Gamma-down}
\end{figure}

As with increasing $t$, the trajectories for decreasing $t$ exhibit
three stages of increasing length. We present the data, much as before, in 
Figures~\ref{short-down}, \ref{long-down}, \ref{3-short-down}, and 
\ref{Gamma-down}. These show 100 trajectories for 5,000 steps,
two trajectories for 1,000,000 steps, close-ups of three 
of the short downwards trajectories, and a histogram of the lengths of 
stage 2, respectively. The long runs show both $\lambda_2$ and the 2-ear count,
while the short runs only show $\lambda_2$. The short runs illustrate what
happens in stages 1 and 2, while the long runs explore stage 3.

However, there are some
important differences between the dynamics for decreasing $t$ and the
dynamics for increasing $t$.
First, because the entropy at $t=t_{low}$ is substantially lower than
the entropy at $t=t_{high}$, the high-to-low transition has to be
forced. As long as $t$ is above the target value, the dynamics only
allow edge swaps that decrease the number of triangles, or at worst
keep the number constant. The initial stage of the transition, rather
than being an entropy-driven relaxation, is essentially a greedy
algorithm for decreasing $t$. This stage is a little slower than the
first stage of the upward transitions, taking around 180 steps instead
of about 50.

\begin{figure}[htbp]
    \includegraphics[width=2in]{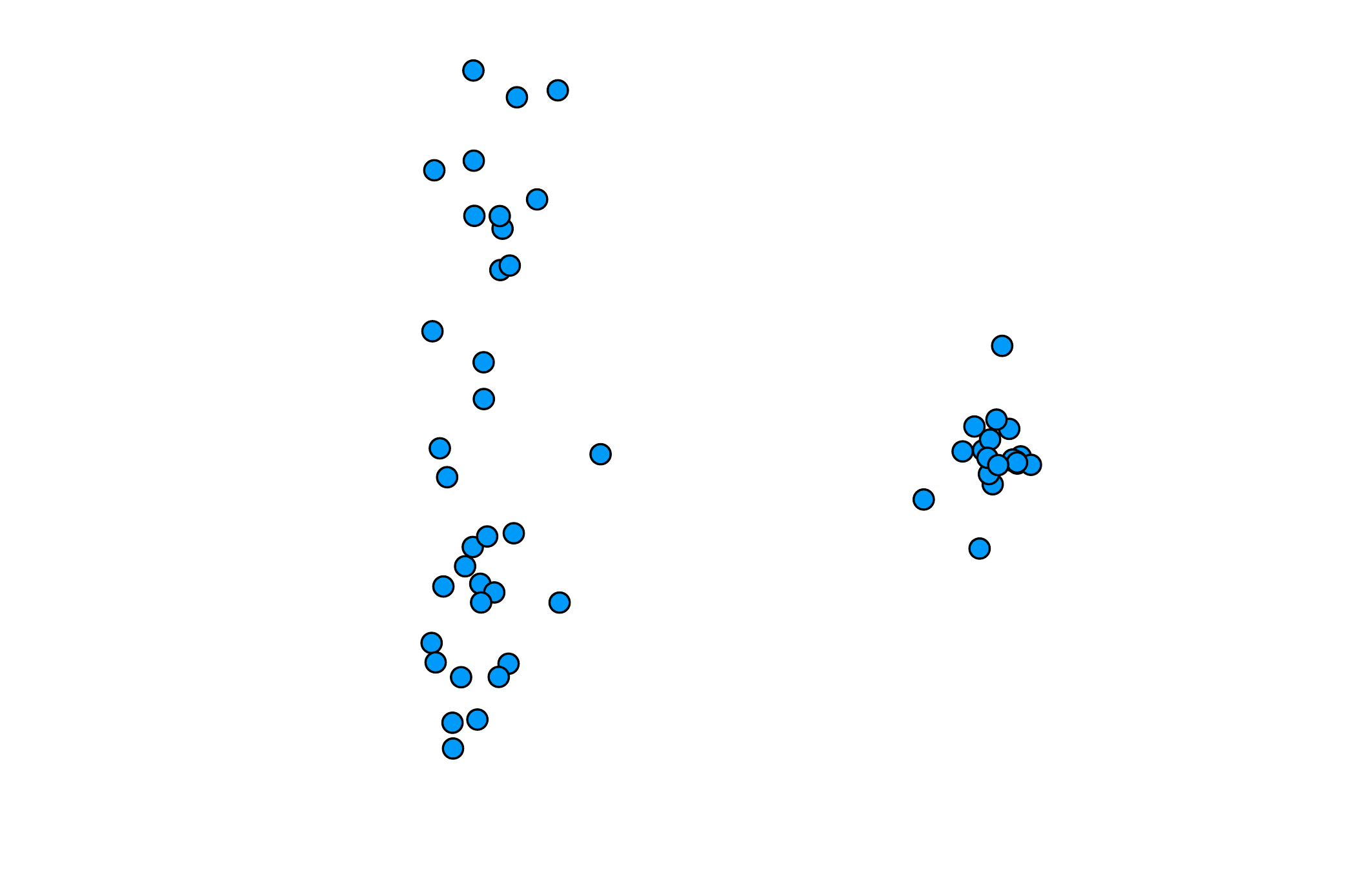}
    \includegraphics[width=2in]{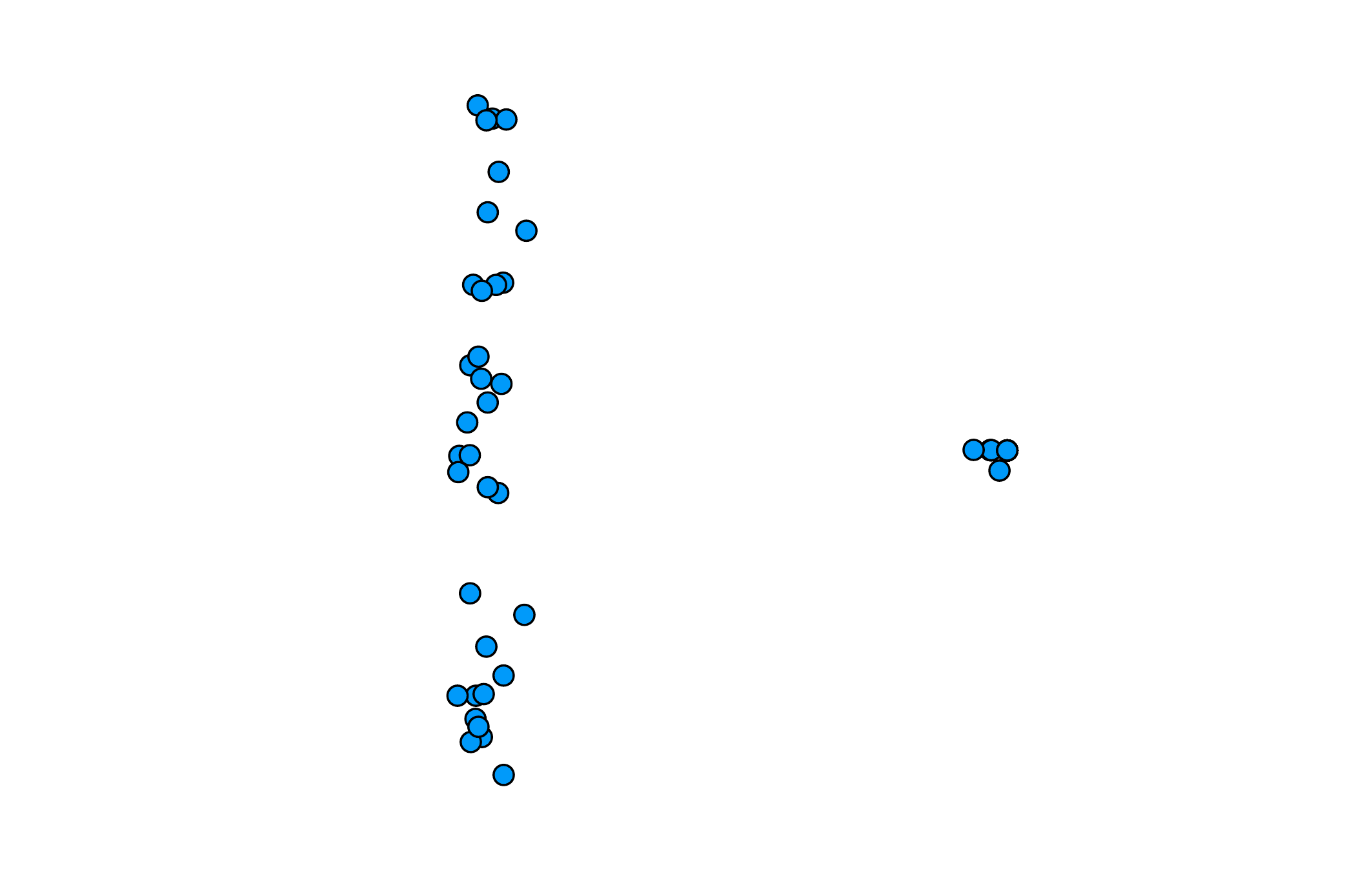}
    \includegraphics[width=2in]{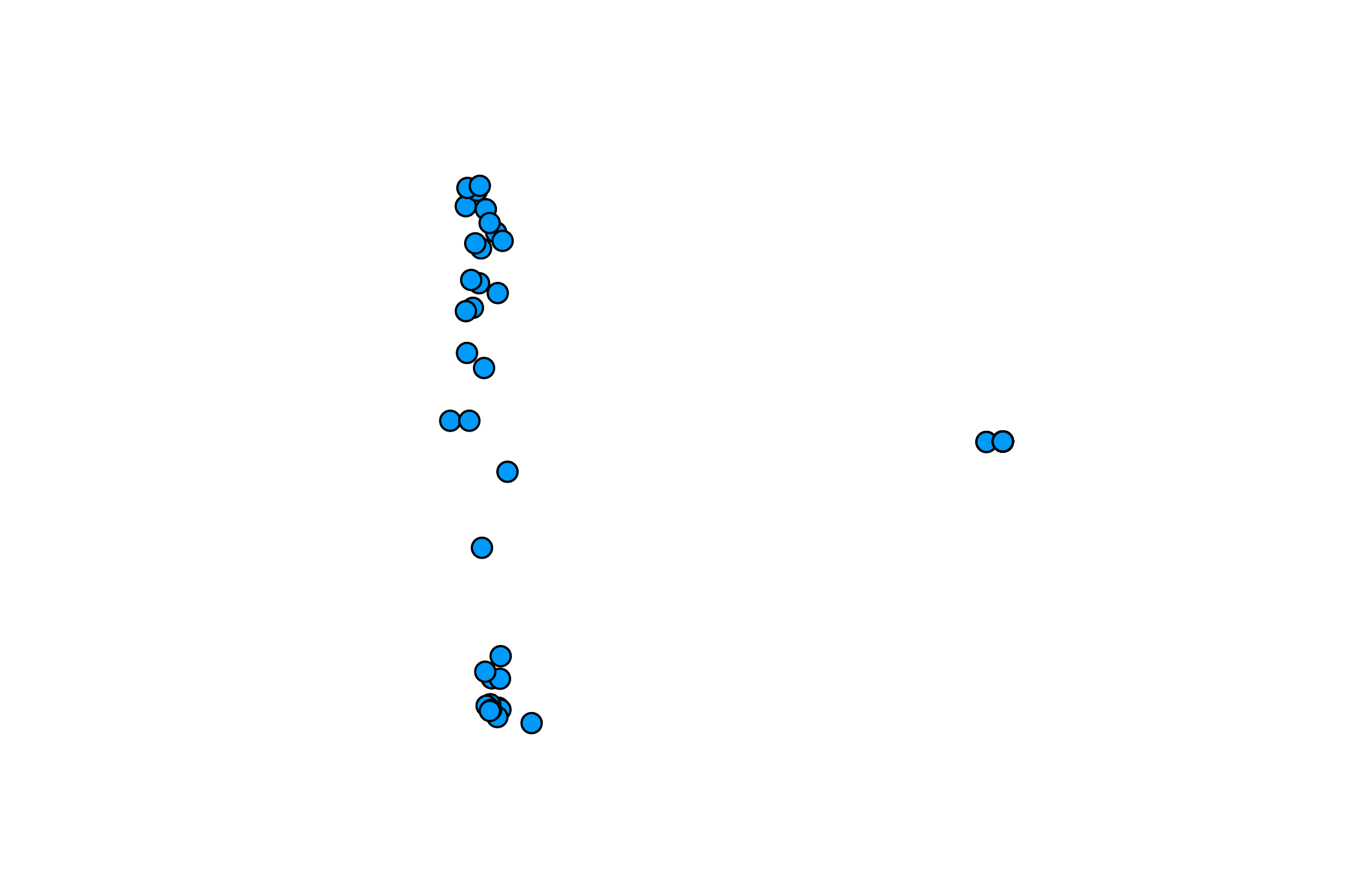}
\caption{Snapshots of $(\xi_1, \xi_2)$ after 10, 100, and 200 steps 
during an downward trajectory. 
Each dot corresponds to a vertex $v$;
the coordinates of the dot are given by $\xi_1(v)$ and $\xi_2(v)$.
The podes are not actually changing size.
Rather, values of $(\xi_1, \xi_2)$  
are becoming so closely positioned that multiple nodes 
only show up as a single dot. This transition can also be seen in the
\href{https://www.youtube.com/playlist?list=PLZcI2rZdDGQqx3WEY6BXoJXcxqQXF5ZzQ}{animations}.}
\label{Xi-down}
\end{figure}

Second, the qualitative form of the graph cannot stay
the same throughout stage 1. Below $t=0.26$, it is mathematically 
impossible to achieve
the desired number of triangles with two podes and with Erd\H{o}s-R\'enyi
structure within each pode. Instead, as noted above, the system has
to develop additional structure within the larger pode. 
As $t$ drops, vertices in this pode
start to segregate into two sub-podes, with edges being more likely
between sub-podes than within a sub-pode. When we reach $t_{low}$,
this segregation is still not complete, but the vertical bar in the $(\xi_1,
\xi_2)$ plot has become more like a dumbbell, with most points near one
end or the other, as can be seen in Figure~\ref{Xi-down}.
At this point, $\lambda_2$ has dropped to around $-13$.

Once we reach $t_{low}$, the dynamics no longer force us to decrease
$t$. The second stage is shorter than for increasing $t$, and 
lasts on the order of 1000 edge swaps. In this stage, 
$\lambda_2$ continues to drop to a value near $-15$ 
as the two sub-podes become better
defined.
(The division of the original
large cluster is not always exactly even, and this affects the resulting
values of $\lambda_2$.)  
However, the decrease in $\lambda_2$ is not as rapid as in stage 1 and is no
longer monotonic. Indeed, the onset of this stage can be viewed not only as the
point where $t$ reaches the target value, but also as the
point where the plot of $\lambda_2$ versus time starts to become noisy. 
This is clearly seen in Figure \ref{3-short-down}, where the 
data has not been smoothed. At the
end of the second stage, our networks now involve two
(more-or-less) identical clusters, each a subset of the original large
pode, and a third cluster that is different, with no evidence of any
internal structure within these three new podes.

There is considerably variation in the length of stage 2, as can be seen
in Figure~\ref{short-down}. As before, we define this length by doing 
a fit of the data after the end of stage 1, and until the 10,000th step, 
by a slanted line followed by a horizontal line. Three such fits (cut off
after the 5000th step) are shown
in Figure~\ref{3-short-down}. In most cases, such as run 9 with a 
long stage 2 and run 46 with a very
short stage 2, the lines fit together clearly, with little
discontinuity. In a few cases, such as run 99, $\lambda_2$ overshoots the
mark and rebounds, or shows other behavior that does not really fit a 
linear model, making it difficult to pinpoint the end of 
stage 2. The distribution of the length of stage 2 is shown in 
Figure~\ref{Gamma-down}, together with a best-fit Gamma functions, with
parameters $\alpha=2.25$ and $\theta=660$. Compared to the similar 
histogram for the upwards trajectory, the value of $\alpha$ is much smaller, 
and the fit is not as good, perhaps because of outliers such as run 99. 

The third stage is much more stable than for increasing $t$. There are
still frequent excursions, where a node attempts to leave its pode and
usually comes back. These appear as short upwards spikes in
$\lambda_2$, with simultaneous downward spikes in the 2-ear
count. Successful migrations are much less common.  For instance, in
the first plot of Figure~\ref{long-down}, the first clear
migration only occurs after about 180,000 steps. There is another one
after 480,000 steps, and a third at about 700,000 steps, plus less
clear events at about 790,000, 830,000, and 960,000 steps. The lowest
values of $\lambda_2$, around $-16.3$, are associated with the most
symmetric configurations, with 18 nodes in each pode, while higher
plateaus around $-15.7$ and $-15.2$ represent 19-18-17 or 20-17-17 or
19-19-16 splits.

\section{Conclusion}\label{conclusion}

In this paper we considered analogues, for constrained random networks, of
the nucleation of phase transitions in molecular systems as normally treated
statistically, in particular the melting/freezing transition between
neighboring fluid and solid phases, which has been difficult to
analyze analytically~\cite{Br,Uh}.

Statistical mechanics produces phase diagrams representing states
which have settled (after `infinite time') into thermal
equilibrium. The states are of macroscopic molecular systems, which
means that one is taking an infinite time limit followed by an
infinite volume limit (to get well-defined transitions), and the
results 
are applicable to
real systems if sufficient care is taken. The `transitions' one sees
in a phase diagram could be understood to represent results taken
through incremental changes of constraints with ample time allowed for
the system to reequilibrate between steps, so-called `quasistatic'
changes. 

At fixed pressure, a fluid material when cooled can remain fluid
far beyond its freezing point unless it is disturbed by complicating
inhomogenities. This phenomenon is studied through `nucleation
theory', which models the manner by which a supercooled fluid
stochastically produces microscale crystalline clusters which, if
large enough, grow to macroscopic crystals. This phenomenon is
understood as resulting from the improbability for a typical
fluid/disordered cluster of molecules to randomly rearrange into a
{\em large enough} (`critical size') crystal cluster.

In this paper we analyzed one of the phase transitions found in
constrained random networks to see if there are similar barriers for a
network in one phase to rearrange structure to another phase. Our
results, detailed in Section~\ref{nucleation}, show distinct
barriers in both directions. We have included our first attempts to analyze these phenomena
statistically.

The current standard for modelling nucleation for molecular systems is
called classical nucleation theory. It is still too crude, with some
predictions wrong by a large factor~\cite{Ca}. It is hard to improve the
theory because one cannot follow nucleation experimentally: it occurs
in very small but unpredictable regions of space, over a very small
time~\cite{Ca}. But versions of `nucleation' can occur in systems
other than molecular matter, and theory can make use of this. Recently
nucleation has been found in physical experiments of cyclically
sheared `sand'~\cite{RRSS,RaSw}, which may prove easier to model
accurately since the system is macroscopic.  Likewise the nucleation
in networks reported here, which was in fact motivated by our work
in~\cite{RRSS,RaSw}, may also contribute to the general
subject. Specifically, as discussed in Section 3 we find that various
stages in the rearrangement of the network structure have well defined
probability distributions, different in the two
directions. Understanding these should be useful for a deeper
nucleation theory.

\end{document}